\documentclass[11pt]{article}
\usepackage{amsmath,enumerate,amsfonts,amssymb,graphicx}
\usepackage{color}
\usepackage{tikz}

\def\RR{{\mathbb R}}
\def\ZZ{{\mathbb Z}}

\def\Sphere{{\mathbb S}}

\newtheorem{theorem} {\sc  Theorem\rm} [section]

\newtheorem{corollary}[theorem] {\sc  Corollary\rm}
\newtheorem{lemma} [theorem] {\sc  Lemma\rm}
\newtheorem{proposition} [theorem] {\sc  Proposition\rm}

\newtheorem{remark}[theorem]{\sc  Remark\rm}

\def\bproof{\noindent{\bf Proof.\;}}
\def\eproof{\hfill$\square$\medskip}

\newcounter{marnote}

\DeclareFontFamily{OT1}{rsfs}{}
\DeclareFontShape{OT1}{rsfs}{m}{n}{ <-7> rsfs5 <7-10> rsfs7 <10-> rsfs10}{}
\DeclareMathAlphabet{\mycal}{OT1}{rsfs}{m}{n}

\def\tr{{\rm tr}}

\def\mcM{{\mycal M}}

\def\mcS{{\mycal S}}
\def\mcC{{\mycal C}}
\def\stg{{\mathfrak{g}}}

\def\tgamma{{\tilde\gamma}}
\def\tils{{\tilde s}}
\def\tilw{{\tilde w}}
\def\ell{{l}}
\def\tpsi{{\tilde \psi}}
\def\stRic{{\overline{\rm Ric}}}
\def\stRiem{{\bar R}}
\def\stGamma{{\bar \Gamma}}
\def\stnabla{{\bar \nabla}}

\newcommand\lr[2]{{\left\langle #1, #2\right\rangle}}

\begin{document}

\title{On smoothness of timelike maximal cylinders in three dimensional vacuum spacetimes}
\author{Luc Nguyen\thanks{Mathematics Department, Princeton University} ~and Gang Tian\thanks{BICMR, Peking University and Mathematics Department, Princeton University}}
\maketitle

\begin{abstract}
We show that timelike maximal cylinders in $\RR^{1 + 2}$ always develop singularities in finite time and that, infinitesimally at a generic singularity, their time slices are evolved by a rigid motion or a self-similar motion. We also prove a mild generalization in non-flat backgrounds.
\end{abstract}

\section{Introduction}

We consider timelike maximal surfaces in a three dimensional vacuum spacetime $(\mcM^{1 + 2},\stg)$. These surfaces are usually referred to as relativistic strings in the literature. 

The case of non-compact timelike maximal graphs in Minkowski spacetimes $\RR^{1 + n}$ is fairly well understood. Global well-posedness for sufficiently small initial data was established by Brendle \cite{Brendle-MaxSurf} and by Lindblad \cite{Lindblad}. The case of general codimension was studied by Allen, Andersson and Isenberg \cite{AllenAI}.

The main focus of this paper is in the case where the surface is an immersed cylinder $\RR \times \Sphere^1$, or a string, in $\RR^{1 + 2}$. We note that local well-posedness for this problem in a larger context was studied recently by Allen, Andersson and Restuccia \cite{AllenAR}.

In \cite{BHNO}, Bellettini, Hoppe, Novaga and Orlandi showed that the problem for a relativistic string in a flat spacetime can be simplified considerably in the so-called orthogonal gauge. In fact, the authors effectively reduced the problem to a homogeneous linear wave equation in one dimension, which admits a simple representation formula from its initial data, which implies in particular the long time existence of parametrizations for relativistic strings. As a consequence, they showed, among other results, that if the initial curve is a centrally symmetric convex curve and the initial velocity is zero, the string shrinks to a point in finite time. (It should be noted that the string does not become extinct there, but rather comes out of the singularity point, evolves back to its original shape and then periodically afterwards.) We also note an earlier paper by Kibble and Turok \cite{KibbleTurok} which showed that a closed string with zero initial velocity must form singularity in finite time.

By a different method, Kong and his collaborators \cite{Kong:2007} and \cite{Kong:2009} proved another representation formula (without the need to fix a gauge). Using their representation formula, they presented many numerical evidences where singularity formation is prominent.

From the above discussion (see also \cite{EggersH}), it is suggestive that, for any arbitrary initial data, a closed string must form singularity in finite time. The main goal of the present paper is to confirm this statement.

Let $\RR^{1 + 2}$ denote the three dimensional Minkowski spacetime endowed with the flat metric, and $(t,x^1, x^2)$ its standard Cartesian coordinates. We prove:

\begin{theorem}\label{Main1}
For any smooth immersed closed curve $\mcC \subset \RR^2 = \{t = 0\}$ and smooth future-directed timelike and nowhere vanishing vector field $V$ along $\mcC$, there exists no globally smooth immersed surface $\mcS \subset \RR^{1 + 2}$ which contains $\mcC$ and tangential to $V$ such that its induced metric is Lorentzian and its mean curvature vector vanishes.
\end{theorem}

Equivalently, the above result asserts that if one evolves a closed curve in $\RR^{1 + 2}$ in a timelike direction such that its mean curvature is zero, it will form singularity in finite time. On the other hand, since the PDE for the parametrization map admits a global solution, it makes sense to talk about the ``maximal surface'' after singularity forms. In Section \ref{Sec:LocPic}, we give a detailed study of the local geometry of a maximal surface at a generic singularity. For a generic singularity propagation, we shows that, locally, the time slices of the maximal surface evolve by a rigid motion: they are either translated or rotated (see Figure \ref{PerSing}). For a generic singularity formation, we show that it is locally self-similar. Self-similar singularity formation was classified by Eggers and Hoppe \cite{EggersH}. Locally, the singularities look like a swallowtail: the first singularity is a cusp of order $4/3$ which splits up to two ordinary cusps at later time (see Figure \ref{SwtailSing}).

As a partial complement to the above theorem, we also establish in Proposition \ref{Prop:LBExistenceTime} a lower bound for the existence time before singularity forms. Our estimate implies, for example, that if $\alpha$ is a non-compact curve in $\RR^2$ such that its total absolute curvature is smaller than one-half, then there is a (possibly immersed) regular timelike maximal surface in $\RR^{1 + 2}$ containing $\alpha$ and perpendicular to $\RR^2$.

In general vacuum spacetimes, the question of how singularities form is less clear. However, using ODE techniques for proving blow-up results of semilinear wave equation (see e.g. \cite{Levine77, John79, Kato80, Glassey81}), we can prove the following result, which implies a singularity statement in $\RR^{1 + 2}$ when the initial curve is convex in $\RR^2$ and the normal velocity vector is parallel along the initial curve.

\begin{theorem}\label{Main2}
Let $(\mcM^{1 + 2},\stg)$ be a complete, oriented, time-oriented, globally hyperbolic, three dimensional vacuum spacetime and $\stnabla$ its connection. Let $\mcC$ be a smooth embedded spacelike acausal closed curve. Along $\mcC$, let $U$ be its unit tangent vector field, $V$ a smooth future-directed timelike unit vector field normal to $U$, and $\nu$ a unit (spacelike) vector field normal to both $U$ and $V$. If
\[
\stg(\stnabla_U U,\nu)^2 - \stg(\stnabla_U V,\nu)^2 > 0 \text{ along } \mcC,
\]
then there exists no globally smooth embedded surface $\mcS^{1+1} \subset \mcM^{1 + 2}$ which contains $\mcC$ and tangential to $V$ such that its induced metric is Lorentzian and its mean curvature vector vanishes.
\end{theorem}

It should be noted that in Theorem \ref{Main2}, we leave out the issue whether the parametrization map exists for all time.

As a final remark, we note that the singularity character of maximal surfaces (as in Theorem \ref{Main1}) is a special feature of the fact that those are surfaces in three dimensions. In the appendix, we give a construction of a regular (two-dimensional) timelike cylindrical maximal surface in $\RR^{1 + 3}$. We suspect that, in $\RR^{1 + n}$ with $n \geq 3$, for generic initial data, timelike maximal cylinders are smooth; but we have not attempted to analyze this statement.

\medskip
{\noindent\bf Acknowledgments.} The authors would like to thank Professor Kong for useful correspondence.

\section{Timelike cylindrical maximal surfaces in $\RR^{1 + 2}$}

The main goal of this section is to prove Theorem \ref{Main1} for timelike maximal surfaces in $\RR^{1 + 2}$.

\subsection{Smooth spatially closed timelike surfaces}\label{ssec:TSurfFlat}

Consider a smooth immersed surface $\mcS$ in $\RR^{1 + 2}$ given by 
\[
\mcS = \Big\{(t = F^0(s^1, s^2), x^1 = F^1(s^1, s^2), x^2 = F^2(s^1, s^2):  (s^1, s^2) \in \omega \subset \RR^2\Big\}
\]
such that the induced metric
\[
g = g_{ab}\,ds^a\,ds^b = \eta(F_{,s^a},F_{,s^b})\,ds^a\,ds^b
\]
is Lorentzian, i.e. $\det g < 0$. Here $a$, $b$ ranges over $\{1,2\}$. Also assume that $\mcC = \mcS \cap \{t = 0\}$ is a smooth immersed closed curve.

We claim that, for any fixed $t$, the cross section
\[
\mcC_t = \mcS \cap \{t = {\rm const}\} 
\]
is either a smooth curve or is empty. To see this, assume that $\mcC_t$ is non-empty and pick $p = (p^1, p^2) \in \mcC_t$. Notice that $(F^0_{,s^1}(p), F^0_{,s^2}(p)) \neq 0$. For if this fails, we must have at $p$ that
\begin{align*}
0 
	&> \det g = \Big[(F^1_{,s^1})^2 + (F^2_{,s^1})^2\Big]\Big[(F^1_{,s^2})^2 + (F^2_{,s^2})^2\Big] - \Big[F^1_{,s^1}\,F^1_{,s_2} + F^2_{,s^1}\,F^2_{,s^2}\Big]^2
		\\
	&= \Big[F^1_{,s^1}\,F^2_{,s_2} - F^1_{,s^2}\,F^2_{,s^1}\Big]^2,
\end{align*}
which is impossible. We thus assume without loss of generality that $F^0_{,s^1}(p) \neq 0$. Then, by the Inverse Function Theorem, we can find a function $f = f(t,s^2)$ such that $p^1 = f(t,p^2)$ and
\[
F^0(f(t,s^2), s^2) = t.
\]
It is thus seen that, near $p$, $\mcC_t$ is given by
\[
\{\gamma^1(t,s^2), \gamma^2(t,s^2)\}
\]
where $\gamma^a(t,s^2) = F^a(f(t,s^2),s^2)$. We will show that this gives a well-parametrized curve. Assume otherwise, then we must have
\[
F^a_{,s^1}\,f_{,s^2} + F^a_{,s^2} = 0.
\]
On the other hand, by definition of $f$,
\[
F^0_{,s^1}\,f_{,s^2} + F^0_{,s^2} = 0.
\]
It follows that
\[
g = \big[-(F^0_{,s^1})^2 + (F^1_{,s^1})^2 + (F^2_{,s^1})^2\big]\Big\{(ds^1)^2 - 2f_{,s^2}\,ds^1\,ds^2 + f_{,s^2}^2\,(ds^2)^2\Big\},
\]
which further implies
\[
\det g = 0,
\]
which violates our assumption that $g$ is Lorentzian. We have thus shown that $\mcC_t$ is a smooth curved.

From the foregoing discussion, $\mcS$ can be represented by
\[
\mcS = \{F(t,s) := (t,\gamma(t,s)) \in \RR^{1 + 2}, T_1 < t < T_2, s \in \RR\}
\]
where $\gamma: \RR_+ \times \RR \rightarrow \RR^2$. Using the timelike character of $\mcS$ is easy to see that $T_1 = - \infty$, $T_2 = + \infty$ and each curve $\mcC_t$ is a closed curve. We thus have
\[
\mcS = \{F(t,s) := (t,\gamma(t,s)) \in \RR^{1 + 2}, t \in \RR, s \in \RR\}
\]
and $\gamma$ is periodic in $s$ with period $\Xi > 0$.

\subsection{The equations}\label{ssec:TheEqns}

We now consider a timelike, topologically cylindrical maximal surface $\mcS$ of the form
\[
\mcS = \{F(t,s) := (t,\gamma(t,s)) \in \RR^{1 + 2}, 0 \leq t < T \leq \infty, s \in \RR\}
\]
and $\gamma$ is periodic in $s$ with period $\Xi > 0$.

The induced metric $g$ is
\begin{equation}
g = -(1 - |\gamma_{,t}|^2)\,dt^2 + 2\,\lr{\gamma_{,t}}{\gamma_{,s}}\,dt\,ds + |\gamma_{,s}|^2\,ds^2,
\label{Dec10-met}
\end{equation}
where $|\cdot|$ and $\lr{\cdot}{\cdot}$ represent the Euclidean norm and dot product of vectors in $\RR^2$. That $g$ is Lorentzian becomes
\begin{equation}
\det g = -|\gamma_{,s}|^2\,(1 - Q) < 0,
\label{Dec10-Lor}
\end{equation}
where
\begin{equation}
Q = |\gamma_{,t}|^2 - \frac{\lr{\gamma_{,t}}{\gamma_{,s}}^2}{|\gamma_{,s}|^2} \geq 0.
\label{Dec10-Q}
\end{equation}
It should be noted that $Q$ remains unchanged under a reparametrization of the form $(t,s) \mapsto  (t,\tilde s(t,s))$.

Let $\nabla$ denote the connection of $g$. The Gauss equation gives
\[
\nabla_{AB} F^\alpha = L_{AB}\,\nu^\alpha
\]
where $L$ is the second fundamental form of $M$ and $\nu$ is the unit normal to $\mcS$. Taking trace with respect to $g$ yields
\[
\Box_g F^\alpha = H\,\nu^\alpha = 0.
\]
Componentwise, we get
\begin{align}
&\Box_g t
	= 0
	\;,\label{Dec10-H1}\\
&\Box_g \gamma
	= 0
	\;.\label{Dec10-H2}
\end{align}

In \cite{BHNO}, Bellettini, Hoppe, Novaga and Orlandi showed that these equations simplify considerably in the so-called orthogonal gauge. For  completeness, we quickly rederive the reduction here.

From the first equation, we have
\begin{equation}
-\partial_t \Big(\frac{|\gamma_{,s}|}{\sqrt{1 - Q}}\Big) + \partial_s \Big(\frac{\lr{\gamma_{,t}}{\gamma_{,s}}}{|\gamma_s|\,\sqrt{1 - Q}}\Big) = 0
	.\label{Dec10-H1*}
\end{equation}
We now look for a reparametrization, say $s = w(t,\tils)$, $\tilde s = \tilw(t,s)$, such that in the new coordinates $(t,\tils)$,
\begin{equation}
\partial_{\tils} \Big(\frac{\lr{\tgamma_{,t}}{\tgamma_{,\tils}}}{|\tgamma_{\tils}|\,\sqrt{1 - Q}}\Big) = 0
	.\label{Dec10-Gauge1}
\end{equation}
In the above $\tgamma(t,\tils) = \gamma(t,s) = \gamma(t,w(t,\tils))$. We have
\[
\tgamma_{,t} = \gamma_{,t} + w_{,t}\gamma_{,s} \text{ and } \tgamma_{,\tils} = w_{,\tils}\,\gamma_{,s}.
\]
As $w(t,\tilw(t,s)) = s$, we also have
\[
w_{,\tils} = \tilw_{,s}^{-1} \text{ and } w_{,t} = - \tilw_{,s}^{-1}\,\tilw_{,t}.
\]
We hence get
\[
\partial_{\tils} \Big(\frac{\lr{\tgamma_{,t}}{\tgamma_{,\tils}}}{|\tgamma_{\tils}|\,\sqrt{1 - Q}}\Big)
	= \tilw_{,s}^{-1}\,\partial_s\Big(\frac{\lr{\gamma_{,t}}{\gamma_{,s}} + \tilw_{,s}^{-1}\,\tilw_{,t}\,|\gamma_{,s}|^2}{|\gamma_{,s}|\,\sqrt{1-Q}}\Big).
\]
From this, it is easy to see that, to achieve \eqref{Dec10-Gauge1}, we solve for $\tilw$ from
\begin{equation}
\left\{\begin{array}{l}
\tilw_{,t} - \mu\,\tilw_{,s} = 0,\\
\tilw(0,s) = s,
\end{array}\right.
	\label{Dec10-Gauge1Cond1}
\end{equation}
where $\mu$ is a solution to
\[
\partial_s\Big(\frac{\lr{\gamma_{,t}}{\gamma_{,s}} + \mu\,|\gamma_{,s}|^2}{|\gamma_{,s}|\,\sqrt{1-Q}}\Big)  = 0.
\]
It is easy to see that
\begin{equation}
\mu = - \frac{\lr{\gamma_{,t}}{\gamma_{,s}}}{|\gamma_{,s}|^2} + \frac{C(t)\,\sqrt{1 - Q}}{|\gamma_{,s}|}
	.\label{Dec10-mu}
\end{equation}

Now \eqref{Dec10-Gauge1Cond1} is a linear transport equation and can be solved easily. The solution $\tilw$ is constant along characteristics which are integral curves of the ODE
\begin{equation}
\dot s(t) = \mu(t,s(t)).
	\label{May09-CharEq}
\end{equation}
Since $\mu$ is smooth and periodic in $s$, $\mu$ is uniformly Lipschitz in $s$ for $t$ belonging to any closed interval of $[0,T)$. It follows that the integral curve of \eqref{May09-CharEq} exists and is smooth and non-crossing for $t \in [0,T)$. It follows that $\tilw$ is strictly increasing in $s$. This shows that $\tils$ is a valid reparametrization. Also, by the periodicity of $\mu$,
\[
\tilw(t,s+\Xi) - \tilw(t,s) = \Xi,
\]
which implies that $\tgamma$ is periodic in $\tils$ with the same period $\Xi$.

In any event, in the new coordinate $(t,\tils)$, \eqref{Dec10-H1*} becomes
\[
\partial_t \Big(\frac{|\tgamma_{,\tils}|}{\sqrt{1 - Q}}\Big) = 0
\]
and so
\begin{equation}
\frac{|\tgamma_{,\tils}|}{\sqrt{1 - Q}} = \rho(\tils).
	\label{Dec10-ConsLaw}
\end{equation}

\begin{proposition}[\cite{BHNO}]\label{BHNO}
Assume that $\mcS$ is a regular timelike maximal surface in $\RR^{1 + 2}$ which is diffeomorphic to $[0, T) \times \Sphere^1$. There exists a smooth parametrization 
\begin{eqnarray*}
[0,T)\times \RR &\rightarrow& \mcS\\
(t,s) &\mapsto& (t,\gamma(t,s))
\end{eqnarray*}
of $\mcS$ such that $\gamma$ is periodic with period $\Xi > 0$ in $s$, $\gamma_{,s}$ is nowhere zero, and
\begin{align}
&\lr{\gamma_{,t}}{\gamma_{,s}} = 0
	,\label{Dec10-Gauge1*}\\
&|\gamma_{,t}|^2 + |\gamma_{,s}|^2 = 1
	,\label{Dec10-ConsLaw*}\\
&\gamma_{,tt} - \gamma_{,ss} = 0
	.\label{Dec11-WaveEq}
\end{align}
Conversely, if $\gamma$ is a regular solution to \eqref{Dec11-WaveEq} and satisfies \eqref{Dec10-Gauge1*} and \eqref{Dec10-ConsLaw*} at initial time then it gives rise to a regular timelike maximal surface in $\RR^{1 + 2}$ for at least some positive time.
\end{proposition}

\bproof We first pick a parametrization of $\mcS$ such that, on the initial curve $\mcC$, $|\gamma_{,s}(0,s)|^2 = 1 - Q(0,s)$. We then define a new coordinate system $(t,\tils)$ by solving \eqref{Dec10-Gauge1Cond1}-\eqref{Dec10-mu} with $C \equiv 0$. Since $\mu$ is smooth, so is the characteristic curves of \eqref{Dec10-Gauge1Cond1}, which implies that $(t,\tils)$ is a smooth coordinate system on $\mcS$ for as long as $(t,s)$ is.

It is straightforward that in the new coordinates, $|\gamma_{,\tils}(0,\tils)|^2 = 1 - Q(0,\tils)$. Equation \eqref{Dec10-Gauge1*}, \eqref{Dec10-ConsLaw*} and \eqref{Dec11-WaveEq} follow from equation \eqref{Dec10-Gauge1}, the conservation law \eqref{Dec10-ConsLaw} and equation \eqref{Dec10-H2}.

For the converse, assume that $\gamma$ satisfies \eqref{Dec11-WaveEq}. We show that if \eqref{Dec10-Gauge1*} and \eqref{Dec10-ConsLaw*} hold at initial time $t = 0$, then they hold everywhere. To this end, consider $\phi = |\gamma_{,t} + \gamma_{,s}|^2$. By \eqref{Dec11-WaveEq},
\[
\phi_{,t}
	= 2\lr{\gamma_{,t} + \gamma_{,s}}{\gamma_{,tt} + \gamma_{,ts}}
	= 2\lr{\gamma_{,t} + \gamma_{,s}}{\gamma_{,ss} + \gamma_{,ts}}
	= \phi_{,s}.
\]
This implies
\[
\phi_{,tt} - \phi_{,ss} = 0.
\]
On the other hand, as \eqref{Dec10-Gauge1*} and \eqref{Dec10-ConsLaw*} hold initially, we have
\[
\phi = 1 \text{ and } \phi_{,t} = \phi_{,s} = 0 \text{ at time } t = 0.
\]
By the uniqueness of the linear wave equation, this implies that $|\gamma_{,t} + \gamma_{,s}|^2 = \phi = 1$ for all time. Similarly, $|\gamma_{,t} - \gamma_{,s}|^2 = 1$ for all time. These two identities imply \eqref{Dec10-Gauge1*} and \eqref{Dec10-ConsLaw*}. The conclusion follows from the discussion preceding the proposition.
\eproof

\subsection{Unavoidable bad parametrization}

We now give the proof of Theorem \ref{Main1}. Assume otherwise that there is a smooth maximal surface $\mcS$ which contains $\mcC$ and tangential to $V$. By Proposition \ref{BHNO}, there is a parametrization $(t,s) \mapsto (t,\gamma(t,s))$ of $\mcS$ with $\gamma$ periodic in $s$, $\gamma_{,s}$ is nowhere zero such that \eqref{Dec10-Gauge1*}, \eqref{Dec10-ConsLaw*} and \eqref{Dec11-WaveEq} hold.

Let $\alpha(s) = \gamma(0,s)$ and $\beta = \gamma_{,t}(0,s)$. Then
\begin{equation}
\lr{\beta}{\alpha_{,s}} = 0 \text{ and } |\alpha_{,s}|^2 + |\beta|^2 = 1 
	.\label{Dec11-Init}
\end{equation}

From the linear wave equation \eqref{Dec11-WaveEq}, we see that
\begin{equation}
\gamma(t,s) = \frac{1}{2}\big(\alpha(s+t) + \alpha(s-t)) + \frac{1}{2}\int_{s-t}^{s + t} \beta(\xi)\,d\xi
	.\label{Dec11-RepF}
\end{equation}

It is readily seen that, if $\alpha(s)$ and $\beta(s)$ are smooth, then $\gamma$ is smooth. Furthermore, if $\alpha$ is a regular parametrization of the initial curve, i.e. $\alpha_{,s}(\cdot)$ is nowhere zero, then for some $T > 0$, $\gamma(t,\cdot)$ gives a regular parametrization for $0 < t < T$. From \eqref{Dec11-RepF}, we see that $\gamma(t,\cdot)$ is not a well parametrization if and only if there exists $s$ such that
\[
\alpha_{,s}(s + t) + \beta(s + t) + \alpha_{,s}(s - t) - \beta(s - t) = 0
	.
\]
We thus introduce
\begin{equation}
a(s) := \alpha_{,s}(s) + \beta(s) \text{ and } b(s) = \alpha_{,s}(s) - \beta(s)
	.\label{Dec27-ab}
\end{equation}

The following lemma shows that bad parametrization always happens, which contradicts the construction of the coordinate system $(t,s)$ of $\mcS$ and thus concludes the proof of Theorem \ref{Main1}.

\begin{lemma}\label{NonemptyBPara}
There exist $s$ and $r$ such that
\begin{equation}
a(s) + b(r) = 0
	.\label{Dec11-SingCond}
\end{equation}
\end{lemma}

\bproof Let $\mcC$ denote the initial curve defined by $\alpha$. Set
\[
A = \Big\{a(s): s\in \mcC\Big\} \text{ and } B = \Big\{-b(s): s\in \mcC\Big\}
	\;.
\]
Evidently, $A$ and $B$ a closed non-empty connected subsets of $\Sphere^1$. If one of them is equal to $\Sphere^1$, \eqref{Dec11-SingCond} holds trivially. Assume thus that $A$ and $B$ are not $\Sphere^1$. Then $A \cup B$ cannot be $\Sphere^1$ as well. Arguing indirectly, assume further that \eqref{Dec11-SingCond} fails. Then there is a connected interval $I \subset \Sphere^1$ such that $A \Subset I$ and $B \Subset \Sphere^1 \setminus I$. This implies that there exist a unit vector $n$ and a real number $\lambda \in (-1,1)$ such that
\[
\lr{p}{n} > \lambda > \lr{q}{n} \text{ for any } p \in A \text{ and } q \in B
	\;.
\]
Therefore, if $L$ is the period of $\alpha$, then
\[
L\,\lambda < \int_0^L \lr{a(s)}{n}\,ds = \int_0^L \frac{d}{ds}\lr{\alpha(s)}{n}\,ds + \int_0^L \lr{\beta(s)}{n}\,ds = \int_0^L \lr{\beta(s)}{n}\,ds
	\;,
\]
and
\[
L\,\lambda > \int_0^L \lr{-b(s)}{n}\,dy = \int_0^L \frac{d}{ds}\lr{-\alpha(s)}{n}\,ds + \int_0^L \lr{\beta(s)}{n}\,ds = \int_0^L \lr{\beta(s)}{n}\,ds
	\;.
\]
This absurdity proves the result.
\eproof

\begin{remark}\label{Com:NEmBPa}
In fact, if $\mcC$ has non-zero rotation index, one can show that, for any $r \in \mcC$, there exists $s = s(r) \in \mcC$ such that \eqref{Dec11-SingCond} holds. To see this, fix $r \in \mcC$. Define
\[
H(\tau,s) = \frac{\alpha_{,s}(s) + \tau\,\beta(s)}{|\alpha_{,s}(s) + \tau\,\beta(s)|}.
\]
By \eqref{Dec11-Init}, $H$ is a continuous map of $[0,1] \times \mcC$ into $\Sphere^1$. Moreover, for $\tau= 0$, $H(0,\cdot)$ is the tangent map of $\mcC$, i.e. it assigns each point of $\mcC$ to the unit tangent vector of $\mcC$ thereof. As $\mcC$ has non-zero rotation index, $H(0,\cdot)$ has non-zero degree. By the homotopy invariance property of the degree, $H(1,\cdot) = a(\cdot)$ also has non-zero degree. In particular, for any $r \in \mcC$, there exists $s \in \mcC$ such that $a(s) = -b(r)$. Here we have used $|b| = 1$.

When $\mcC$ has zero rotation index, it is impossible to have a strong conclusion as in the previous paragraph. For example, when $\beta \equiv 0$, there exists $r_0 \in \mcC$ such that $a(s) + b(r_0) = \alpha_{,s}(s) + \alpha_{,s}(r_0) \neq 0$ for all $s \in \mcC$.
\end{remark}

\subsection{General analysis of the (spatial) unit tangent map}

In the rest of the section, we analyze the local picture where the orthonormal gauge fails. To analyze what happens when the parametrization goes badly, we look into the unit tangent map
\[
U(t,s) = \frac{\gamma_{,s}(t,s)}{|\gamma_{,s}(t,s)|}.
\]
In terms of $a$ and $b$,
\begin{equation}
U(t,s) = \frac{a(s + t) + b(s - t)}{|a(s + t) + b(s - t)|^{3}}
	.\label{Dec26-UTangM}
\end{equation}

We note that \eqref{Dec11-Init} implies that $|a| = |b| = 1$. Thus, there exist smooth functions $\zeta$ and $\eta$ such that
\begin{equation}
a_{,s}(s) = \zeta(s)\,a^\perp(s) \text{ and } b_{,s}(s) = \eta(s)\,b^\perp(s)
	.\label{Dec27-xieta}
\end{equation}
where the perpendicular rotation $v^\perp$ of a vector $v = (v_1, v_2)$ is defined as $v^\perp = (-v_2, v_1)$.

\begin{lemma}\label{BDCritW}
If $\gamma_{,s}(t_0,s_0) = 0$, i.e.
\begin{equation}
a(s_0 + t_0) + b(s_0 - t_0) = 0
	,\label{Dec27-SingCondXX}
\end{equation}
and if
\begin{equation}
\zeta(s_0 + t_0) \neq \eta(s_0 - t_0)
	\label{Dec30-SingCondKey}
\end{equation}
then the unit tangent map $U(t_0,\cdot)$ of the curve $\mcC_{t_0}$ is discontinuous at $s_0$. More specifically, it reverses direction across $s_0$.
\end{lemma}

\bproof Let $e = a(s_0 + t_0) = - b(s_0 + t_0)$, $s_0^+ = s_0 + t_0$ and $s_0^- = s_0 - t_0$. By \eqref{Dec27-xieta},
\begin{align*}
a(s + t_0) + b(s - t_0)
	&= (s - s_0)\,[\zeta(s_0^+) - \eta(s_0^-)]e^\perp
			+ O(|s - s_0|^2)
\end{align*}
where the big $O$ notation is meant for $s$ close to $s_0$. This implies that
\begin{align}
U(t_0,s) 
	&= \frac{(s - s_0)\,[\zeta(s_0^+) - \eta(s_0^-)]}{|s - s_0||\zeta(s_0^+) - \eta(s_0^-)|}\,e^\perp
		 + O(|s - s_0|).	
	\label{Jan05-X1}
\end{align}
This shows that $U(t_0, s)$ reverses direction as $s$ changes across $s_0$.
\eproof

\begin{remark}\label{rem:genericsing}
In fact, under the hypotheses of Lemma \ref{BDCritW},
\begin{align*}
\gamma(t_0,s) 
	&= \gamma(t_0,s_0) + \frac{1}{2}\,(s - s_0)^2\,[\zeta(s_0^+) - \eta(s_0^-)]\,e^\perp\\
		&\qquad\qquad + \frac{1}{6}\,(s - s_0)^3 \Big\{[\zeta_{,s}(s_0^+) - \eta_{,s}(s_0^-)]e^\perp
				- [\zeta^2(s_0^+) - \eta^2(s_0^-)]\,e\Big\}.
\\		&\qquad\qquad + \frac{1}{24}(s - s_0)^4\,\Big\{[\zeta_{,ss}(s_0^+) - \eta_{,ss}(s_0^-)
			- \zeta^3(s_0^+) + \eta^3(s_0^-)]e^\perp\\
			&\qquad\qquad\qquad\qquad - 3[\zeta(s_0^+)\,\zeta_{,s}(s_0^+)
				- \eta(s_0^-)\,\eta_{,s}(s_0^-)]\,e\Big\}
\end{align*}
Thus, if $\zeta^2(s_0^+) \neq \eta^2(s_0^-)$, the curve $\mcC_{t_0}$ has an ordinary cusp at $s_0$.
\end{remark}

\begin{remark}\label{rem:swallowtail}
We claim that if $t_0 > 0$ is the smallest time such that $\gamma_{,s}(t_0, s_0) = 0$ for some $s_0$, then
\[
\zeta(s_0 + t_0) = \eta(s_0 - t_0).
\]
Assume this claim for the moment and assume in addition that
\[
\zeta_{,s}(s_0 + t_0) \neq \eta_{,s}(s_0 - t_0).
\]
Then the curve $\mcC_{t_0}$ has a cusp of order $4/3$ at $s_0$. Furthermore, for $t > t_0$, this singularity splits up into two ordinary cusps. This picture is consistent with \cite{EggersH}. See Section \ref{Sec:LocPic} for a more detailed discussion.
\end{remark}

To prove the claim above, note that
\[
\partial_s |a(s + t) + b(s - t)|^2 = \lr{a^\perp(s + t)}{b(s - t)}(\zeta(s + t) - \eta(s - t)).
\]
and
\[
\partial_s^2\Big|_{(t,s) = (t_0,s_0)} |a(s + t) + b(s - t)|^2 = (\zeta(s_0^+) - \eta(s_0^-))^2.
\]
Thus, by the implicit function theorem, if $\zeta(s_0^+) \neq \eta(s_0^-)$ then there exist some $\epsilon > 0$ and a smooth map $S: (t_0 - \epsilon, t_0 + \epsilon) \rightarrow \RR$ such that $S(t_0) = s_0$ and 
\[
0 = \partial_s |a(S(t) + t) + b(S(t) - t)|^2 = \lr{a^\perp(S(t) + t)}{b(S(t) - t)}(\zeta(S(t) + t) - \eta(S(t) - t)).
\]
Also, as $\zeta(s_0^+) \neq \eta(s_0^-)$, we can also assume that $\zeta(S(t) + t) \neq \eta(S(t) - t)$ for $t \in (t_0 - \epsilon, t_0 + \epsilon)$. This implies that
\[
\lr{a^\perp(S(t) + t)}{b(S(t) - t)} = 0 \text{ for } t \in (t_0 - \epsilon, t_0 + \epsilon).
\]
As $a(s_0^+) + b(s_0^-) = 0$, the continuity of $a$ and $b$ implies that $2\gamma_{,s}(t,S(t)) = a(S(t) + t) + b(S(t) - t) = 0$, which contradicts our assumption on $t_0$.

\subsection{The case of non-zero rotation index}

In the following discussion, we write
\[
a = (\cos\psi, \sin \psi) \text{ and } b = -(\cos \tpsi, \sin\tpsi).
\]
Note that $\psi' = \zeta$ and $\tpsi' = \eta$.

Since $a$ and $b$ are periodic and having the same degree, say $d$, (see Remark \ref{Com:NEmBPa}), we have
\[
\psi(s + L) - \psi(s) = \tpsi(r + L) - \tpsi(r) = 2\,d\,\pi
\]
where $L$ is the period of $\alpha$. Also, since 
\[
\alpha_{,s} = \frac{1}{2}( a + b) = \sin\frac{\psi - \tpsi}{2} \Big(- \sin \frac{\psi + \tpsi}{2}, \cos \frac{\psi + \tpsi}{2}\Big)
\]
and $\alpha_{,s}$ is nowhere vanishing, we infer that the range of $\psi - \tpsi$ does not intersect $2\pi\ZZ$. 

\begin{lemma}\label{NZI::DoubleImFlat}
Assume that $\alpha$ has non-zero rotation index and the unit tangent map is continuous. If $\psi(s_0) = \psi(s_1) = \tpsi(r_0)$ for some $s_0$, $s_1$ and $r_0$ with $0 < s_1 - s_0 < L$, then $\psi$ is constant in $(s_0, s_1)$.
\end{lemma}

\bproof We will only consider the case where the rotation index $d$ of $\alpha$ is positive. The other case can be proved similarly.

Arguing indirectly, we assume that $\psi$ is non-constant in $(s_0, s_1)$. Then 
\begin{equation}
\text{either }\max_{[s_0,s_1]} \psi > \psi(s_0) \text{ or } \min_{[s_0,s_1]} \psi < \psi(s_0).
	\label{Jan15-1}
\end{equation}

Assume for now that the former case holds. Set
\[
M = \min\Big(2d\pi,\max_{[s_0,s_1]} \psi - \psi_0\Big) > 0,
\]
and define
\[
s_- = \sup\{ s < s_1: \psi(s) = \psi(s_0) + M/4\} \text{ and } s_+ = \inf\{s > s_-: \psi(s) = \psi(s_0) + 3M/4\}.
\]
By the mean value theorem, there exists $s_2 \in [s_-,s_+]$ such that
\[
\psi'(s_2) = \frac{\psi(s_+) - \psi(s_-)}{s_+ - s_-} > 0.
\]
By definition of $s_\pm$, we also have $\psi(s_0) + M/4 < \psi(s_2) \leq \psi(s_0) + 3M/4$. Since $\tpsi(r_0) = \psi(s_0)$ and $\tpsi(r_0 + L) = \psi(s_0) + 2d\pi$, the intermediate value theorem implies that there exists $r_2 \in [r_0, r_0 + L]$ such that $\tpsi(r_2) = \psi(s_2)$. Now let
\[
s_3 = \sup\{s: \psi(\hat s) > \psi(s_2) \text{ for all } s_2 < \hat s < s\}.
\]
Since $\psi'(s_2) > 0$ and $\psi(s_1) = \psi(s_0) < \psi(s_2)$, $s_3$ exists and $s_2 < s_3 < s_1$. Furthermore, $\psi(s_3) = \psi(s_2)$ and $\psi'(s_3) \leq 0$.

We thus end up with
\[
a(s_2) = a(s_3) = -b(r_2), \zeta(s_2) > 0 \geq \zeta(s_3).
\]
Therefore,
\[
\text{either } \zeta(s_2) \neq \eta(r_2) \text{ or } \zeta(s_3) \neq \eta(r_2).
\]
Then Lemma \ref{BDCritW} applies yielding that $\kappa$ must blow up somewhere, a contradiction.

If the second case in \eqref{Jan15-1} holds, the argument is similar using the comparison values in of $\tpsi(r)$ for $r \in [r_0 - L, r_0]$.
\eproof 

\begin{corollary}\label{NZI::Mon}
Assume that $\alpha$ has non-zero rotation index and $\kappa$ is always finite. Then $\psi$ and $\tpsi$ are either both non-increasing or non-decreasing.
\end{corollary}

\bproof Again, we will only consider the case where the rotation index $d$ of $\alpha$ is positive. By Remark \ref{Com:NEmBPa}, there exists $r_0$ such that $a(0) + b(r_0) = 0$. We can further assume that $\tpsi(r_0) = \psi(0)$.

We claim that $\psi[0,L] \subset [\psi(0),\psi(L)]$. Define
\[
s_+ = \inf\{s > 0: \psi(s) = \psi(L)\} \leq L \text{ and } s_- = \sup\{s < L: \psi(s) = \psi(0)\} \geq 0.
\]
Since $\psi(0) = \psi(s_-) = \tpsi(0)$, Lemma \ref{NZI::DoubleImFlat} shows that $\psi$ is constant in $(0,s_-)$. Similarly, $\psi$ is constant in $(s_+,L)$. The claim follows easily.

We now show that $\psi$ is non-decreasing. Assume otherwise, then for some $0 \leq s_0 < s_1 \leq L$, $\psi(s_0) > \psi(s_1)$. By the claim, $s_0 > 0$. Thus, by the intermediate value theorem, there exists $s_2 \in (0,s_0)$ such that $\psi(s_2) = \psi(s_1)$. By Lemma \ref{NZI::DoubleImFlat}, $\psi$ is constant in $(s_2,s_1)$ contradicting the assumption that $\psi(s_0) > \psi(s_1)$. We hence conclude that $\psi$ is non-decreasing. Similarly, $\tpsi$ is non-decreasing.
\eproof

\begin{proposition}\label{NZI::NoSmooth}
For any smooth initial data $\alpha$ and $\beta$ such that $\alpha$ has non-zero rotation index, there exists a time $T$ such that the curvature $\kappa(T,\cdot)$ of the curve $\mcC_T$ blows up.
\end{proposition}

\bproof We will only consider the case where the rotation index $d$ of $\alpha$ is positive.

Assume for some initial data $\alpha$ and $\beta$ that the curvature function $\kappa$ of the solution $\gamma$ remains finite for all time. This implies in particular that the unit tangent map $U(t,s)$ is a continuous function. We will show that the curve $\gamma(t,\cdot)$ will contract to a point in finite time, which results in a contradiction.

By \eqref{Dec11-RepF},
\[
U(t,s) = {\rm sgn}\Big(\sin \frac{\psi(s + t) - \tpsi(s - t)}{2}\Big) \Big(- \sin \frac{\psi(s + t) + \tpsi(s - t)}{2}, \cos\frac{\psi(s + t) + \tpsi(s - t)}{2}\Big),
\]
where $\rm sgn$ denotes the sign function. It follows that, for any $t$, 
\begin{equation}
\text{the function $\psi(\cdot + t) - \tpsi(\cdot - t)$ does not change sign.}
	\label{Jan15-2}
\end{equation}
For otherwise, $U$ must be discontinuous there.

By Corollary \ref{NZI::Mon}, $\psi$ and $\tpsi$ are both non-decreasing. (Here we have also used the fact that $\psi(L) - \psi(0) = 2\pi\,d > 0$.)

Using Remark \ref{Com:NEmBPa}, the mean and intermediate value theorems as in the proof of Lemma \ref{NZI::DoubleImFlat}, we can find $s_0$ and $t_0$ such that
\begin{equation}
\psi(s_0 + t_0) = \tpsi(s_0 - t_0),
	\label{Jan15-3}
\end{equation}
and
\begin{equation}
\psi'(s_0 + t_0) = \tpsi'(s_0 - t_0) > 0.
	\label{Jan15-4}
\end{equation}

By \eqref{Jan15-2}, we have either
\begin{equation}
\psi(s + t_0) \geq \tpsi(s - t_0) \text{ for all } s
	\label{Jan15-5}
\end{equation}
or
\begin{equation}
\psi(s + t_0) \leq \tpsi(s - t_0) \text{ for all } s.
	\label{Jan15-5*}
\end{equation}

Assume for now that \eqref{Jan15-5} holds. By \eqref{Jan15-3} and \eqref{Jan15-4}, for some $\delta_0 > 0$, there holds
\[
\psi(s_0 + t_0 - \delta) < \tpsi(s_0 - t_0 + \delta) \text{ for all } \delta \in (0,\delta_0).
\]
Thus, by \eqref{Jan15-2}
\begin{equation}
\psi(s + t_0 - \delta) \leq \tpsi(s - t_0 + \delta) \text{ for all } s \text{ and for all } \delta \in (0,\delta_0).
	\label{Jan15-6}
\end{equation}
From \eqref{Jan15-5} and \eqref{Jan15-6} we deduce that
\[
\psi(s + t_0 - \delta) \leq \tpsi(s - t_0 + \delta) \leq \psi(s + t_0 + \delta) \text{ for all } s \text{ and for all } \delta \in (0,\delta_0).
\]
Sending $\delta \rightarrow 0$, we thus get
\[
\psi(s + t_0) \equiv \tpsi(s - t_0).
\]
This implies that 
\[
\gamma_{,s}(t_0,s) = a(s + t_0) + b(s - t_0) \equiv 0,
\]
which shows that $\mcC_{t_0}$ is actually a point.

The case where \eqref{Jan15-5*} holds can be handled similarly. We get
\[
\psi(s + t_0 + \delta) \geq \tpsi(s - t_0 - \delta) \geq \psi(s + t_0 - \delta) \text{ for all } s \text{ and for all } \delta \in (0,\delta_0).
\]
This again forces $\psi(s + t_0) \equiv \tpsi(s - t_0)$ and thereby concludes the proof.
\eproof

\subsection{The case of zero rotation index}

We next switch to the case where the rotation index of $\alpha$ is zero. We have
\begin{lemma}\label{ZI::NoPlateau}
Assume that $\alpha$ has zero rotation index and $U$ is continuous. For any $r$, there is no more than one $s$ such that $\psi(s) = \tpsi(r)$.
\end{lemma}

\bproof Assume by contradiction that there exists $s_0$, $s_1$ and $r_0$ with $0 < s_1 - s_0 < L$ such that $\psi(s_0) = \psi(s_1) = \tpsi(r_0)$. Define
\[
M = \max_{[s_0, s_0+L]} \psi \text{ and } m = \min_{[s_0, s_0 + L]} \psi.
\]
Since $\psi - \tpsi$ is nowhere zero, $\psi$ is not a constant function. Thus, either $M > \psi(s_0)$ or $m < \psi(s_0)$. In the sequel, we will assume that $M > \psi(s_0)$. The case where $m < \psi(s_0)$ can be treated similarly. Furthermore, by replacing $(s_0, s_1)$ by $(s_1, s_0 + L)$ if necessary, we can assume that either $M$ is achieved in $[s_0, s_1]$.

We claim that $\tpsi(r) \leq \psi(s_0)$ for all $r$. If this is wrong, we can argue using the mean and intermediate value theorems as in the proof of Lemma \ref{NZI::DoubleImFlat} to find $s_2$, $s_3$ and $r_2$ such that
\[
\psi(s_2) = \psi(s_3) = \tpsi(r_2) \text{ and } \zeta(s_2) > 0 \geq \zeta(s_3).
\]
This gives a violation to the conclusion of Lemma \ref{BDCritW}. The claim follows.

Now, consider the interval $(s_1, s_0 + L)$. If the minimum value of $\psi$ in this interval is less than $\psi(s_0)$, the same argument leads to another violation of Lemma \ref{BDCritW}. We thus have
\begin{equation}
\psi(s) \geq \psi(s_0) \geq \tpsi(r) \text{ for any $s$ and $r$}.
	\label{Jan15-X1}
\end{equation}

We next show that
\begin{equation}
\psi(s) - \tpsi(r) \leq 2\pi \text{ for any $s$ and $r$}.
	\label{Jan15-X2}
\end{equation}
Arguing indirectly, assume that \eqref{Jan15-X2} fails. By the intermediate value theorem, there exists $s_2$ and $r_2$ such that $\psi(s_2) = \tpsi(r_2) + 2\pi$. Furthermore, we can assume that $s_2 \in (s_0, s_1)$. Evidently, if $\psi(s_2) = \psi(s_0)$, we can further use the intermediate value theorem again to find $s_2'$ and $r_2'$ such that $\psi(s_0) < \psi(s_2') = \tpsi(r_2') + 2\pi$. We thus assume that $\psi(s_2) > \psi(s_0)$. If $\psi(s_2) < M$, the intermediate value theorem implies that there exists $s_3 \in (s_2, s_1)$ such that $\psi(s_3) = \psi(s_2)$. The argument leading to \eqref{Jan15-X1} then implies that $\psi$ can only takes value either on $(-\infty, \psi(s_2)]$ or $[\psi(s_2),\infty)$, which is obviously not the case. We thus get
\[
\psi(s) = M \text{ whenever there exists $r$ such that } \psi(s) = \tpsi(r) + 2\pi.
\]
Since $\psi$ achieves values in $(\psi(s_0),M)$, this implies that $\tpsi(r) \geq M - 2\pi$, which implies \eqref{Jan15-X2}, a contradiction. We have thus shown \eqref{Jan15-X2}.

We now revisit the proof of Lemma \ref{NonemptyBPara}. Define $A = \{a(s)\}$ and $B = \{-b(r)\}$. By \eqref{Jan15-X1} and \eqref{Jan15-X2}, $A$ and $B$ intersects at exactly two points: 
\[
A \cap B = \Big\{(\cos \psi(s_0),\sin\psi(s_0)), (\cos M, \sin M)\Big\}.
\]
This implies that there exist a unit vector $n$ and a real number $\lambda \in (-1,1)$ such that
\[
\lr{p}{n} \geq \lambda \geq \lr{q}{n} \text{ for any } p \in A \text{ and } q \in B
	\;.
\]
Furthermore, since neither $\psi$ nor $\tpsi$ are constant, there exists $s$ and $r$ such that
\[
\lr{a(s)}{n} > \lambda > \lr{-b(s)}{n}
	\;.
\]
We can then argue as in the proof of Lemma \ref{NonemptyBPara} to reach a contradiction.
\eproof

\begin{proposition}\label{ZI::NoSmooth}
For any smooth initial data $\alpha$ and $\beta$ with $\alpha$ has zero rotation index, there exists a time $T$ such that the unit tangent map $U(T,\cdot)$ of the curve $\mcC_T$ is discontinuous.
\end{proposition}

\bproof Assume that the curvature function remains finite for all time.

By Lemma \ref{NonemptyBPara}, there exists $s_0$ and $r_0$ such that $a(s_0) + b(r_0) = 0$. We can thus assume that $\psi(s_0) = \tpsi(r_0)$. 

Let
\[
M = \max_{[s_0, s_0+L]} \psi \text{ and } m = \min_{[s_0, s_0 + L]} \psi.
\]

We claim that $\psi(s_0) \in \{M, m\}$. Assume otherwise that $m < \psi(s_0) < M$. Fix $0 < s_+ - s_- < L$ such that $\psi(s_-) = m$ and $\psi(s_+) = M$. Then, by the intermediate value theorem, there exists $s_1 \in (s_-, s_+)$ and $s_2 \in (s_+, s_- + L)$ such that $\psi(s_1) = \psi(s_2) = \psi(s_0) = \psi(r_0)$. This violates the conclusion of Lemma \ref{ZI::NoPlateau}. The claim follows.

In the sequel, we assume that $\psi(s_0) = m$. The other can be handled similarly.

By symmetry, $\tpsi(r_0)$ is also an extremal value of $\tpsi$. If it is the minimal value, we can find $s_1$ and $r_1$ such that $\psi(s_1) = \tpsi(r_1)$ and the common value is not extremal, which is a contradiction to the above claim. Thus, 
\[
\tpsi(r_0) = \max \tpsi,
\]
which implies
\[
\psi(s) \geq \tpsi(r) \text{ for any $s$ and $r$}.
\]

As in the proof of Lemma \ref{ZI::NoPlateau}, we next show that
\[
\psi(s) - \tpsi(r)  \leq 2\pi \text{ for any $s$ and $r$}.
\]
If this was not correct, we can find $s_1$ and $r_1$ such that $\psi(s_1) = \hat \psi(r_1)$ where $\hat\psi = \tpsi + 2\pi$. Using the intermediate value theorem, we can further assume that $\psi(s_1) > m$. Again, $\psi(s_1)$ is an extremal value of $\psi$, which must be the maximal value. Similarly, $\hat\psi(r_1)$ is the minimal value of $\hat\psi$. We thus get
\[
0 = \max \psi - \min \hat\psi = \max \psi - \min \tpsi - 2\pi.
\]

We can now argue as in the proof of Lemma \ref{ZI::NoPlateau} to get a contradiction.
\eproof

\subsection{A lower bound for the blow up time}

Having proved a singularity statement, we would like to see how long a solution stays smooth before it develops singularity. The estimate should depends on the initial curve, $\alpha_*: [p,q] \rightarrow \RR^2$ and the initial (normal) velocity field $\partial_t + \beta_*$ along $\alpha_*$. To clarify the notation, $\alpha = \gamma(0,s)$ and $\beta(0,s) = \gamma_{,t}(0,s)$ are reparametrizations of $\alpha_*$ and $\beta_*$, i.e. $\alpha = \alpha_* \circ \Phi$ and $\beta = \beta_* \circ \Phi$ for some diffeomorphism $\Phi$. Furthermore, the estimate should be local, because the speed of propagation is finite for the wave equation. For this latter point, in this section, we do not assume that $\alpha_*$ is a closed curve.

It is useful to define the timelikeness index of the (prospective) maximal surface along $\alpha_*$ to be
\[
j(\alpha_*,\mcS) = j(\alpha_*,\beta_*) := L(\alpha_*)\left\{\int_{\alpha_*} \frac{1}{\sqrt{-|V_*|^2}}\,|d\alpha_*|\right\}^{-1} = L(\alpha_*)\left\{\int_{\alpha_*} \frac{1}{\sqrt{1 - |\beta_*|^2}}\,|d\alpha_*|\right\}^{-1}
\]
where $L(\alpha_*)$ denotes the length of $\alpha_*$. Note that by definition $\beta_*$ is normal to $\alpha_*$ and has norm smaller than $1$. Hence
\[
j(\alpha_*,\beta_*) \in (0,1].
\]

\begin{proposition}\label{Prop:LBExistenceTime}
Let $\alpha_*: [p,q] \rightarrow \RR^2$ be a (not necessarily closed) smooth curve in $\RR^2 = \{t = 0\} \subset \RR^{1 + 2}$, $U_*$ its unit tangent vector field, and $V_* = \partial_t + \beta_*$ a smooth timelike vector field along $\alpha_*$ and normal to $\alpha_*$. Let $\partial_t + a_*$ and $-\partial_t + b_*$ be the null vector fields belonging to the span of $\{U_*,V_*\}$ such that
\[
\lr{a_*}{U_*} > 0 \text{ and } \lr{b_*}{U_*} > 0.
\]
If the timelikeness index along $\alpha_*$ and the curvatures of $a_*$ and $b_*$ satisfy
\[
j(\alpha_*,\beta_*) > \frac{3}{2}\int_{\alpha_*} \left[|\nabla_{U_*} a_*| + |\nabla_{U_*} b_*| \right]\,|d\alpha_*|,
\]
then there exist two smooth functions $p, q: [0,T] \rightarrow \RR$ with $T = L(\alpha_*)/j(\alpha_*,\beta_*)$, $p(0) = p$, $q(0) = q$, $p(T) = q(T)$ and a map $\gamma: \Omega \rightarrow \RR^2$ with $\Omega = \{(t,x): t \in [0,T], x \in [p(t),q(t)]\}$ such that the map $(t,x) \mapsto (t,\gamma(t,x))$ defines a regular timelike maximal surface which contains $\alpha_*$, is tangential to $V_* = \partial_t + \beta_*$ and whose lateral boundary are two null curves.
\end{proposition}

\bproof Switching to isothermal gauge as before, we can drop the subscript $*$. Note that $a_*$ and $b_*$ coincide with the vector field $a$ and $b$ defined in \eqref{Dec27-ab}.

We will show that the map
\begin{align*}
\Omega = \{(t,s): t \in [0,T], p + t \leq s \leq q - t\} \rightarrow \RR^{2 + 1}
(t,s) \mapsto (t,\gamma(t,s))
\end{align*}
defines a smooth maximal surface. (Note that $2T = |q - p|$ in this gauge.) To this end, it suffices to show that $\gamma_{,s}(t,s) \neq 0$ for $(t,s) \in \Omega$. 

We first estimate $|\alpha_{,s}| = |\gamma_{,s}(0,\cdot)|$. Using the function $\xi$ and $\eta$ defined in \eqref{Dec27-xieta}, we estimate for $x,y \in [p,q]$:
\begin{align*}
|\alpha_{,s}(x)| - |\alpha_{,s}(y)| 
	&= \int_x^y \frac{\lr{\alpha_{,s}(z)}{\alpha_{,ss}(z)}}{|\alpha_{,s}(z)|}\,dz\\
	&= \int_x^y \frac{\lr{a(z) + b(z)}{a_{,s}(z) + b_{,s}(z)}}{2|a(z) + b(z)|}\,dz\\
	&= \int_x^y \frac{\lr{a(z) + b(z)}{a^\perp(z)\,\xi(z) + b^\perp(z)\,\eta(z)}}{2|a(z) + b(z)|}\,dz\\
	&= \int_x^y [\xi(z)-\eta(z)]\frac{\lr{a^\perp(z)}{b(z)}}{2|a(z) + b(z)|}\,dz\\
	&= \frac{1}{4}\int_x^y [\xi(z)-\eta(z)]\lr{\frac{a(z) + b(z)}{|a(z) + b(z)|}}{\underbrace{a^\perp(z) - b^\perp(z)}_{=\beta^\perp(z)}}\,dz.
\end{align*}
This implies that
\[
\Big||\alpha_{,s}(x)| - |\alpha_{,s}(y)| \Big| 
	\leq \frac{1}{4}\int_p^q |\xi(z) - \eta(z)|\,|\beta(z)|\,dz
	.
\]
Integrating in the $y$ variables, it follows that
\[
\Big| |\alpha_{,s}(x)| - \underbrace{\frac{L}{2T}}_{= \frac{1}{2}j}\Big| = \Big| |\alpha_{,s}(x)| - \frac{1}{2T}\int_p^q |\alpha_{,s}(y)|\,dy\Big|
	\leq \frac{1}{4}\int_p^q |\xi(z) - \eta(z)|\,|\beta(z)|\,dz.
\]
We thus deduce that
\begin{equation}
|\alpha_{,s}(x)| \geq \frac{1}{2}j - \frac{1}{4}\int_p^q |\xi(z) - \eta(z)|\,|\beta(z)|\,dz \text{ for } x \in [p,q].
	\label{11Oct11-alphasBnd}
\end{equation}

Next, for any $(t,s) \in \Omega$ we have
\begin{align*}
|\alpha_{,s}(t,s)| 
	&= \frac{1}{2}|a(s + t) + b(s - t)|\\
	&\geq \frac{1}{2}|\alpha_{,s}(s - t)| - \frac{1}{2}|a(s + t) - a(s - t)|\\
	&\geq \frac{1}{2}|\alpha_{,s}(s - t)| - \frac{1}{2}\left|\int_{s - t}^{s + t} a^\perp(z)\,\xi(z)\,dz\right|\\
	&\geq \frac{1}{2}|\alpha_{,s}(s - t)| - \frac{1}{2}\int_p^q |\xi(z)|\,dz.
\end{align*}
By symmetry, we have
\[
|\alpha_{,s}(t,s)| 
	\geq \frac{1}{2}|\alpha_{,s}(s - t)| - \frac{1}{2}\int_p^q |\eta(z)|\,dz,
\]
which implies
\[
|\alpha_{,s}(t,s)| 
	\geq \frac{1}{2}|\alpha_{,s}(s - t)| - \frac{1}{4}\int_p^q \,dz.
\]
Combining with \eqref{11Oct11-alphasBnd}, we obtain
\[
|\alpha_{,s}(t,s)| \geq \frac{1}{4}j - \frac{1}{4}\int_p^q [|\xi(z)| + |\eta(z)|][\frac{1}{2}|\beta(z)| + 1]\,dz \text{ for } (t,s) \in \Omega
\]
Now notice that
\[
\nabla_U a = \frac{1}{|\alpha_{,s}|} a_{,s} = \frac{1}{|\alpha_{,s}|} a^\perp\,\xi \text{ and } \nabla_U b = \frac{1}{|\alpha_{,s}|} b^\perp\,\eta,
\]
we arrive at
\[
|\alpha_{,s}(t,s)| \geq \frac{1}{4}j - \frac{1}{4}\int_{\alpha} [|\nabla_U a(z)| + |\nabla_U b(z)|][\frac{1}{2}|\beta(z)| + 1]\,|d\alpha| \text{ for } (t,s) \in \Omega.
\]
Note that as $\partial_t + \beta$ is timelike, $|\beta| < 1$. Hence, by hypothesis, the right hand side of the above inequality is positive. We conclude the proof.
\eproof

As a consequence of the above result we have
\begin{corollary}
For any closed curve $\mcC \subset \RR^2 = \{t = 0\}\subset \RR^{1 + 2}$ and any future-directed timelike field $V$ along $\mcC$, there exist a constant $T_* > 0$ and a regular timelike maximal surface $\mcS$ containing $\mcC$ and tangential to $V$ in the time slab $\{0 \leq t < T_*\}$. Furthermore, the maximal existence time $T_*$ is finite and, for any positive $l$ which is smaller than the length of the final curve $\gamma(T_*,\cdot)$, there holds
\[
\lim_{t \rightarrow T} \sup\Bigg\{ \frac{1}{j(\Gamma,\mcS)}\,\int_\Gamma \left[|\nabla_{U} a| + |\nabla_{U} b| \right]\,|d\Gamma|\Bigg\} > \frac{2}{3}.
\]
where the supremum is taken over all connected sub-arc $\Gamma$ of length $l$ of the curve $\gamma(t,\cdot)$, where $U$ is the unit tangent to $\Gamma$, $\partial_t + a$ and $-\partial_t + b$ are the null vector field along $\Gamma$ which is tangential to $\Gamma$. In particular, the surface $\mcS$ becomes null somewhere on the final curve.
\end{corollary}

\section{Local picture at a singularity in $\RR^{1 + 2}$}\label{Sec:LocPic}

In this section, we study the local picture at a singularity. Let $\alpha: [-1,1] \rightarrow \RR^2$ and $\beta: [-1,1] \rightarrow \RR^2$ be two smooth map such that
\begin{enumerate}[(a)]
  \item $\alpha$ defines a continuous curve which is smooth away from $\alpha(0) = 0$,
  \item $\lr{\alpha'}{\beta} = 0$ and $|\alpha'|^2 + |\beta|^2 = 1$ in $[-1,1]$,
  \item $|\beta(0)| = 1$.
\end{enumerate}
Note that (c) implies 
\begin{equation}
\lr{\beta'(0)}{\beta(0)} = 0 \text{ and } \lr{\beta''(0)}{\beta(0)} + |\beta'(0)|^2 \leq 0,
	\label{Jan12-X1}
\end{equation}
and (b) and (c) imply
\begin{equation}
\alpha'(0) = 0 \text{ and } \lr{\alpha''(0)}{\beta(0)} = 0.
	\label{Jan12-X2}
\end{equation}
In particular, as $\beta(0) \neq 0$,
\[
\text{$\alpha''(0)$ and $\beta'(0)$ are colinear.}
\]
In addition, (b) implies that
\begin{align}
&|\alpha''(0)|^2 + \lr{\beta''(0)}{\beta(0)} + |\beta'(0)|^2 = 0,
	\label{Jan12-X3}\\
&\lr{\alpha'''(0)}{\beta(0)} + 2\lr{\alpha''(0)}{\beta'(0)} = 0.
	\label{Jan12-X4}
\end{align}

In view of \eqref{Dec11-RepF}, define
\[
\gamma(t,s) = \frac{1}{2}(\alpha(s + t) + \alpha(s - t)) + \frac{1}{2}\int_{s - t}^{s + t} \beta(\xi)\,d\xi
\]
for 
\[
(t,s) \in \Omega := \big\{ (t,s): |t| + |s| \leq 1\big\}. 
\]
As shown earlier, $\gamma$ defines a regular timelike maximal surface $\mcS$ away from points where $\gamma_{,s}(t,s) = 0$. We would like to analyze its local behavior near $\gamma(0,0)$.

We note that a spacetime dilation of a maximal surface remains a maximal surface. Thus, it would be natural to consider the limit $n\gamma(\frac{t}{n},\frac{s}{n})$ as $n \rightarrow \infty$. (Note that the rescaling of the parametrization variables is to ensure the gauge conditions \eqref{Dec10-Gauge1*} and \eqref{Dec10-ConsLaw}.) It is easy to see that the limit is the map $(t,s) \mapsto \beta(0)t$ which parametrizes a null plane in $\RR^{1 + 2}$. This approximation of the original maximal surface is rather crude. In what to follow, we would like to obtain a better description.

We start with an analysis of the zero set of $\gamma_{,s}$. We have
\[
\gamma_{,s}(t,s) = \frac{1}{2}(\alpha'(s + t) + \alpha'(s - t)) + \frac{1}{2}(\beta(s + t) - \beta(s - t)) = \frac{1}{2}(a(s + t) + b(s - t)).
\]
where $a$ and $b$ are defined in \eqref{Dec27-ab}. Recalling the function $\zeta$ and $\eta$ defined in \eqref{Dec27-xieta}, we have
\begin{align*}
4\partial_s |\gamma_{,s}|^2(t,s)
	&= \lr{a'(s + t) + b'(s - t)}{a(s + t) + b(s - t)}\\
	&= \lr{a^\perp(s + t)}{b(s - t)}(\zeta(s + t) - \eta(s - t)).
\end{align*}
and
\[
4 \partial_{s}^2 |\gamma_{,s}|^2(0,0) = (\zeta(0) - \eta(0))^2 = 4|\alpha''(0)|^2.
\]
%

\subsection{Generic singularity propagation}\label{ssec:LPicPer}

Let us first consider the case $\alpha''(0) \neq 0$. Note that this implies in particular that the curve defined by $\alpha$ has a cusp at the $\alpha(0)$. In this case, we can find some $\epsilon > 0$ such that the solutions to $\partial_s |\gamma_{,s}|^2(t,s) = 0$ in $(-\epsilon,\epsilon)^2 \subset \Omega$ is given by some smooth curve $\Gamma = \{(t,S(t)): t \in (-\epsilon,\epsilon)\}$. Furthermore, as $\alpha''(0) \neq 0$, $\zeta(S(t) + t) - \eta(S(t) - t) \neq 0$ in $(-\epsilon,\epsilon)^2$, and so $\lr{a^\perp(S(t) + t)}{b(S(t) - t)} = 0$. Continuity then implies that $2\partial_{,s} \gamma(S(t),t) = a(S(t) + t) + b(S(t) - t) = 0$. In this case, locally around $\gamma(0,0)$, the singularities of $\mcS$ are given by $\{(t, \gamma(t,S(t)))\}$, which is null and tangent to $\partial_t + \beta(0)$. Also, as 
\[
4 \partial_t\partial_{s} |\gamma_{,s}|^2(0,0) = \zeta(0)^2 - \eta(0)^2 = 4\lr{\alpha''(0)}{\beta'(0)},
\]
we also have
\[
S'(0) = -\frac{\lr{\alpha''(0)}{\beta'(0)}}{|\alpha''(0)|^2}.
\]
We have thus shown:

\begin{proposition}
Let $\mcS$ be a ``timelike maximal'' surface defined by a smooth map $(t,s) \mapsto \gamma(t,s)$ satisfying $\lr{\gamma_{,t}}{\gamma_{,s}} = 0$ and $|\gamma_{,t}|^2 + |\gamma_{,s}|^2 = 1$. If $\gamma(t_0,s_0)$ is a singular point of $\mcS$ (i.e. $\gamma_{,s}(t_0, s_0) = 0$) and if $\gamma_{,ss}(t_0,s_0) \neq 0$, then locally around $\gamma(t_0,s_0)$ singularities of $\mcS$ are cusps and propagate along the null curve $t \mapsto (t, \gamma(t,S(t))$ where $S$ solves
\[
\left\{\begin{array}{l}
S'(t) = -\frac{\lr{\gamma_{,ss}(t,S(t))}{\gamma_{,ts}(t,S(t))}}{|\gamma_{,ss}(t,S(t))|^2},\\
S(t_0) = s_0.
\end{array}\right.
\]
\end{proposition}

\begin{figure}[h]
\begin{center}
\includegraphics[width=.4\textwidth]{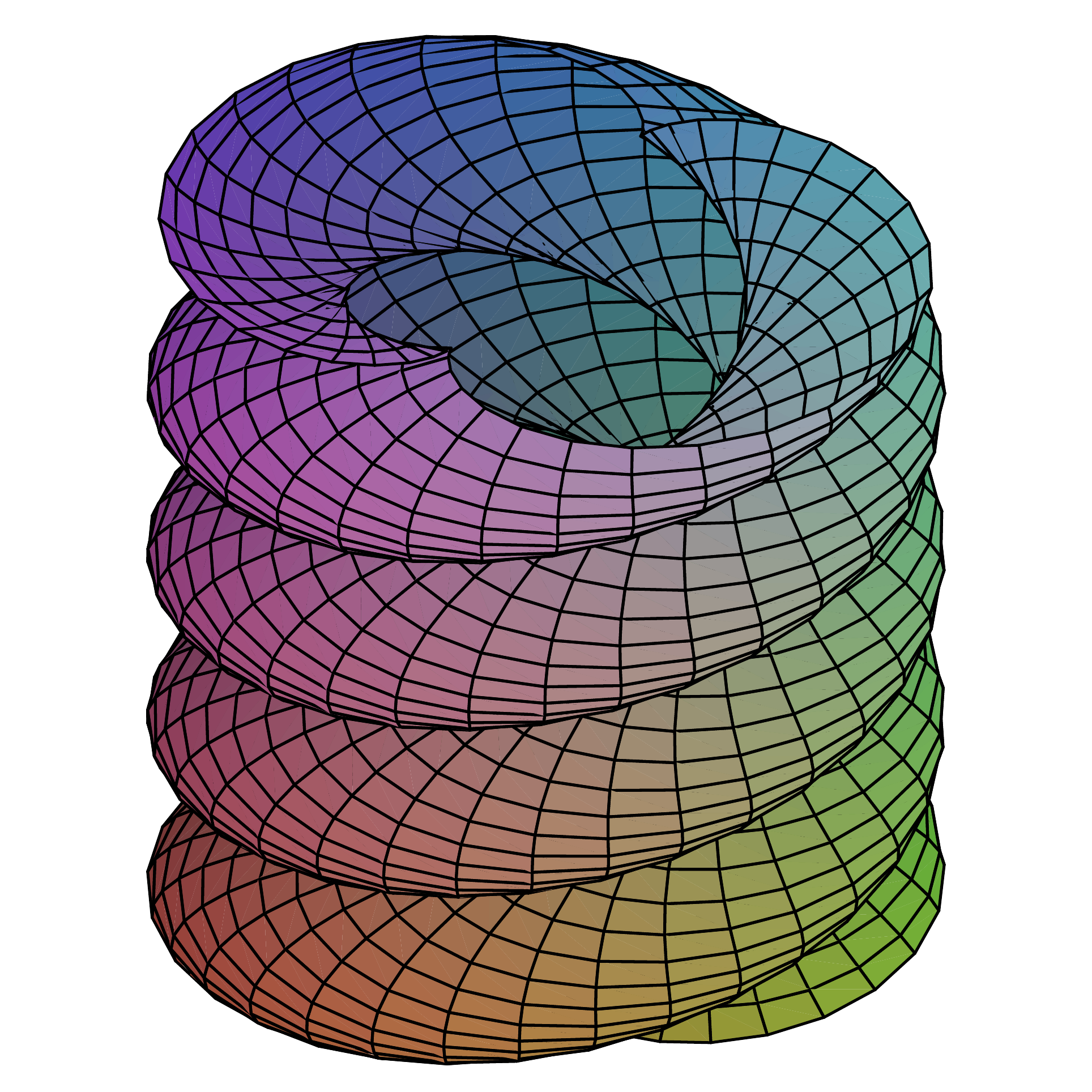}
\caption{A maximal surface in which any time slice has exactly two singularities which travel along two null curves.}
\label{PerSing}
\end{center}
\end{figure}

Examples of exact solutions with the above behavior are given by
\[
\gamma(t,s) = \frac{1}{2}\Big(\frac{1}{\lambda_1} + \frac{1}{\lambda_2} - \frac{\cos \lambda_1(s - t)}{\lambda_1} - \frac{\cos \lambda_2(s + t)}{\lambda_2}, - \frac{\sin \lambda_1(s - t)}{\lambda_1} + \frac{\sin \lambda_2(s + t)}{\lambda_2}\Big),
\]
and ones obtained by sending $\lambda_1 \rightarrow 0$ or $\lambda_2 \rightarrow 0$. The picture of the corresponding $\mcS$ for $\lambda_1 = 3$ and $\lambda_2 = 1$ is given in Figure \ref{PerSing}. In this example, every time slice has two cusp singularities. Those singularities propagate along null helices. As an evolution of curves in $\RR^2$, it is a rotation: the time slice at time $t$ is obtained by rotating the initial curve by $3t$ radian around some point.

Let us show that the behavior seen above is prominent: Infinitesimally around a singularity of the present type, the time slices $\mcC_t$ of $\mcS$ are evolved by a rigid motion, namely either a translation or a rotation. The key idea is that the curve $t \mapsto \gamma(t,S(t))$ of zeroes of $\gamma_{,s}$ can be approximated up to second order around $t = 0$ by its osculating circle.

We note that, up to cubic terms, the Taylor expansion of $\gamma$ around $(0,0)$ is
\[
\gamma(t,s) = \beta(0)\,t + \frac{1}{2}\alpha''(0)(t^2 + s^2) + \beta'(0)\,ts + \frac{1}{6}\,\beta''(0)\,(t^3 + 3\,ts^2) + \frac{1}{6}\,\alpha'''(0)\,(3t^2 s + s^3) + \ldots
\]
Let
\[
e = \beta(0), p = \lr{\alpha''(0)}{e^\perp} \text{ and } q = \lr{\beta'(0)}{e^\perp}.
\]
Using \eqref{Jan12-X1}-\eqref{Jan12-X4}, we have
\begin{align*}
\gamma(t,s) 
	&= t\,e - \frac{1}{6}\,[(p^2 + q^2)(t^3 + 3\,t\,s^2) + 2pq(3t^2\,s + s^3)]e\\
		&\qquad + \frac{1}{2}[p(t^2 + s^2) + 2q\,ts]\,e^\perp + O(|t|^4 + |s|^4) e + O(|t|^3 + |s|^3)e^\perp.
\end{align*}
In particular,
\begin{align*}
\alpha(s)
	&= -\frac{pq}{3}\,s^3\,e + \frac{p}{2}\,s^2\,e^\perp + O(s^4)\,e + O(s^3)\,e^\perp.
\end{align*}

The curvature of the curve $t \mapsto \gamma(t,S(t))$ at $t = 0$ is
\[
k_0 = \lr{\alpha''(0) - \frac{\lr{\alpha''(0)}{\beta'(0)}}{|\alpha''(0)|^2}\,\beta'(0)}{\beta(0)^\perp} = p - \frac{q^2}{p}.
\]

If $k_0 \neq 0$, i.e. $p \neq \pm q$, then
\begin{align*}
\gamma(t,s) 
	&= \frac{1}{2}\Big(-\frac{1}{p-q}\,\sin [(p-q)(s - t)] + \frac{1}{p + q}\sin [(p + q)(s + t)]\Big)\,e\\
		&\qquad + \frac{1}{2}\Big(\frac{2p}{p^2 - q^2} - \frac{1}{p - q}\,\cos [(p-q)(s - t)] - \frac{1}{p + q}\cos [(p + q)(s + t)] \Big)e^\perp\\
		&\qquad + O(|t|^4 + |s|^4) e + O(|t|^3 + |s|^3)e^\perp.
\end{align*}
Hence, if we set $X(t,s) = \lr{\gamma(t,s)}{e^\perp}  - \frac{p}{p^2 - q^2}$ and $Y(t,s) = \lr{\gamma(t,s)}{e}$, we have
\begin{align*}
X(t,s)
	&= - \sin k_0 t\,X\big(0,s  + \frac{q}{p} t\big) + \cos k_0 t\,Y\big(0,s + \frac{q}{p}t\big)
		+ O(|t|^3 + |s|^3)e^\perp,\\
Y(t,s)
	&= \cos k_0 t\,X\big(0,s  + \frac{q}{p} t\big) + \sin k_0 t\,Y\big(0,s + \frac{q}{p}t\big)
		+ O(|t|^4 + |s|^4) e.
\end{align*}
This shows that, infinitesimally, $\mcC_t$ ``is'' the image of a rotation of $\mcC_0$ by $\frac{q}{p}t$ radian about $\frac{p}{p^2 - q^2}e^\perp$ (which is the center of the osculating circle of the curve of zeroes of $\gamma_{,s}$ at $s = t = 0$).

The case where $k_0 = 0$ can be obtain by considering the limit $p \rightarrow q$ or $p \rightarrow - q$. For example, when $p = q$, we have
\begin{align*}
\gamma(t,s) 
	&= t\,e + \frac{1}{3}\,p^2(s + t)^3e
		 + \frac{1}{2}\,p\,(s + t)^2\,e^\perp + O(|t|^4 + |s|^4) e + O(|t|^3 + |s|^3)e^\perp\\
	&= t\,e + \alpha(s \pm t) + O(|t|^4 + |s|^4) e + O(|t|^3 + |s|^3)e^\perp.
\end{align*}
This shows that, infinitesimally, the curve $\mcC_t$ ``is'' a translation of $\mcC_0$. The exact solution approximant is
\[
\gamma^*(t,s) = \frac{1}{2}\big(t - s + \frac{1}{2p}\,\sin 2p(s + t)\big)e + \frac{1}{4p}\big(1 - \cos 2p(s + t)\big)e^\perp.
\]
%

\subsection{Generic singularity formation}\label{ssec:LPicSw}

Let us consider next the case $\alpha''(0) = 0$. Note that condition (b) together with $\alpha''(0) = 0$ implies that
\[
\text{$\alpha'''(0)$ and $\beta'(0)$ are colinear.}
\]

To make the situation not too degenerate, we make an empirical ansatz that
\begin{equation}
\beta'(0) \neq 0 \text{ and } \alpha'''(0) \neq 0.
\label{5Jan12-EmA}
\end{equation}
As we have said earlier, this is the generic case for singularity formation.

Arguing as in the previous case but considering zero of $\partial_t |\gamma_{,s}|^2$ instead, we see that the singularities of $\mcS$ around $\gamma(0,0)$ are given by a curve $\Gamma =  \{(T(s),s)\}$ where $T$ is smooth and
\[
T'(s) = -\frac{\lr{\alpha''(s)}{\beta'(s)}}{|\beta'(s)|^2}.
\]
Furthermore, note that
\begin{equation}
T'(0) = 0 \text{ and } T''(0) = -\frac{\lr{\alpha'''(0)}{\beta'(0)}}{|\beta'(0)|^2} \neq 0,
	\label{23Jan12-J1}
\end{equation}
which implies that the singularities lie either all in the future or in the past. In addition, by Remark \ref{rem:swallowtail}, $\gamma(t_0,s_0)$ is a cusp of order $4/3$ and other singularities are regular cusps. We thus have:

\begin{proposition}
Let $\mcS$ be a ``timelike maximal'' surface defined by a smooth map $(t,s) \mapsto \gamma(t,s)$ satisfying $\lr{\gamma_{,t}}{\gamma_{,s}} = 0$ and $|\gamma_{,t}|^2 + |\gamma_{,s}|^2 = 1$. If $\gamma(t_0,s_0)$ is a singular point of $\mcS$ (i.e. $\gamma_{,s}(t_0, s_0) = 0$) and if $\gamma_{,ts}(t_0,s_0) \neq 0$, 
then locally around $\gamma(t_0,s_0)$ singularities of $\mcS$ lie along the null curve $s \mapsto (T(s),\gamma(T(s),s))$ where $T$ solves
\[
\left\{\begin{array}{l}
T'(s) = \frac{\lr{\gamma_{,ss}(T(t),s)}{\gamma_{,ts}(T(s),s)}}{|\gamma_{,ts}(T(t),s)|^2},\\
T(s_0) = t_0.
\end{array}\right.
\]
Furthermore, if $\gamma_{,ss}(t_0,s_0) = 0$ and $\gamma_{,sss}(t_0,s_0) \neq 0$, then those singularities are cusps (except possibly $\gamma(t_0,s_0)$) and the curve $s \mapsto (T(s),\gamma(T(s),s))$ lies either all in the past or in the future (depending on whether $\lr{\gamma_{,sss}}{\gamma_{,t}}(t_0,s_0)$ is positive or negative, respectively) and splits into two null curves which join together at $(t_0,\gamma(t_0,s_0))$ as a cusp.
\end{proposition}

\begin{figure}[h]
\begin{center}
\begin{tabular}{|c|c|}
\hline
\includegraphics[width=.4\textwidth]{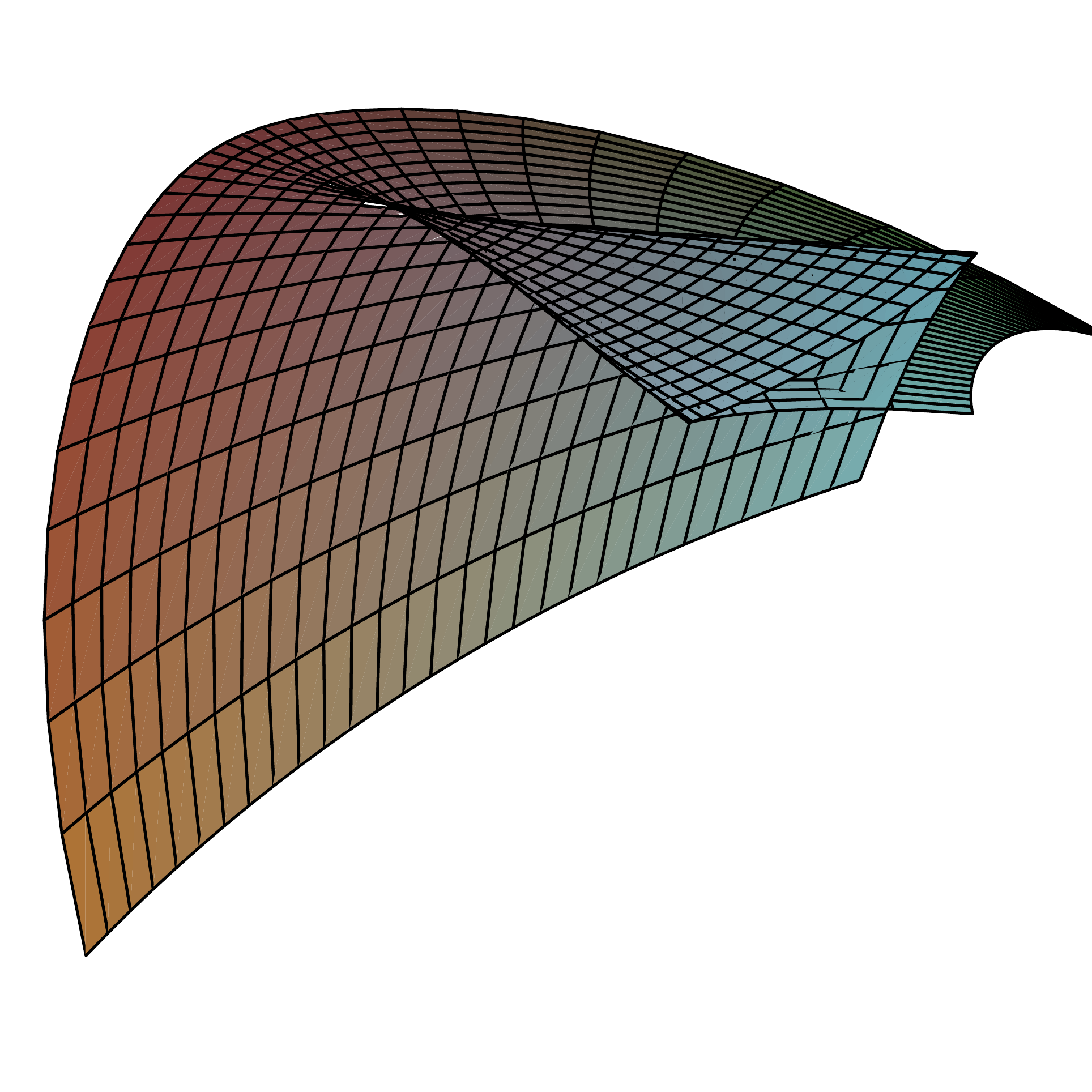}&\includegraphics[width=.4\textwidth]{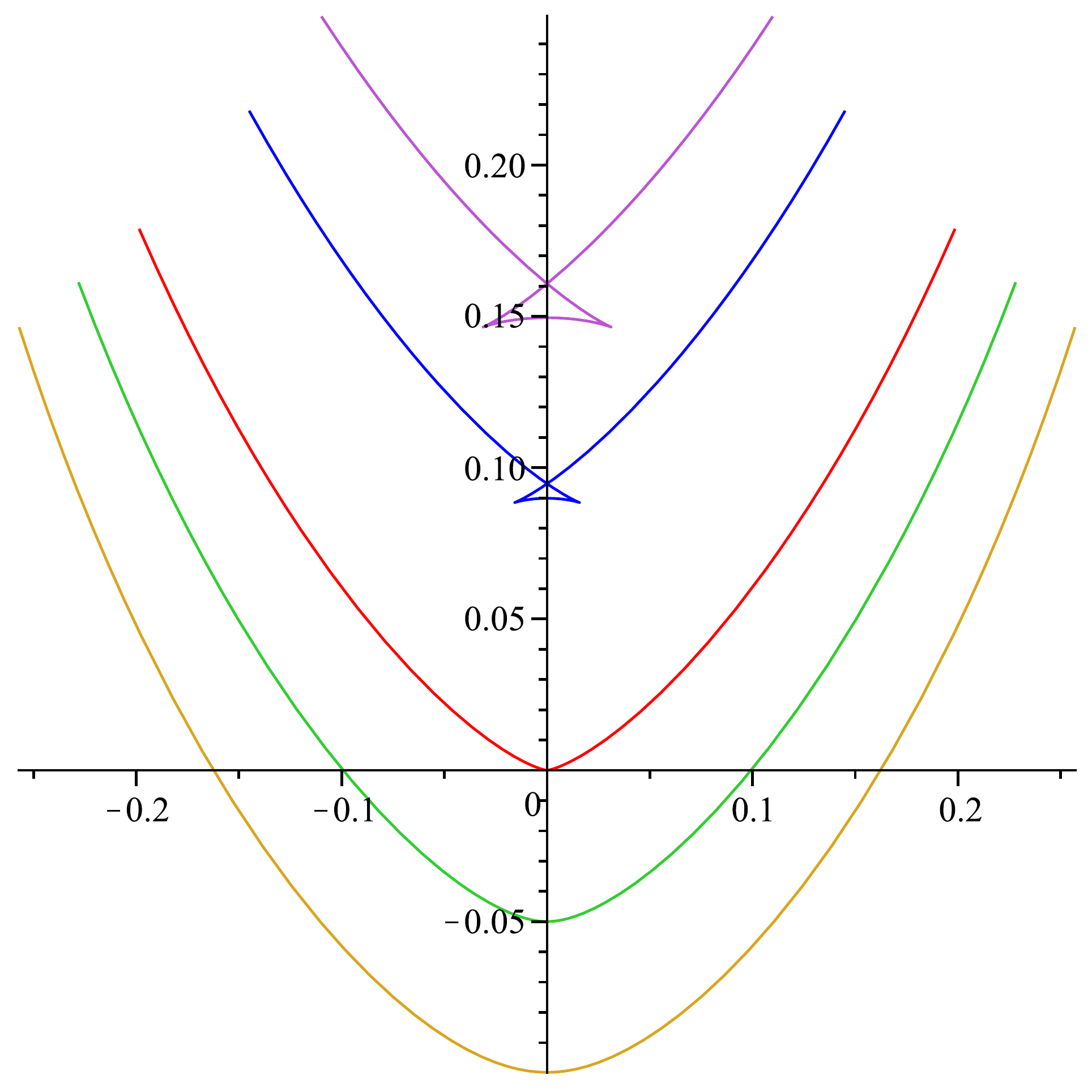}\\
(a) & (b)\\
\hline
\end{tabular}
\caption{(a) A maximal surface with swallowtail-type singularity. (b) Time slices of the maximal surface in (a).}
\label{SwtailSing}
\end{center}
\end{figure}

An example is given by
\begin{align*}
\alpha(s) 
	&= \Big(\frac{1}{3}\,\sin^3 s, \frac{2}{3} -\frac{2}{3}\cos s + \frac{1}{3}\sin^2 s\,\cos s\Big),\\
\beta(s)
	&= (-\sqrt{1 - \sin^4 s}\sin s, \sqrt{1 - \sin^4 s}\,\cos s).
\end{align*}
The picture of the corresponding $\mcS$ is given in Figure \ref{SwtailSing}.

As we discussed before, the set of singularities looks like a swallowtail. Under a self-similar assumption, this was proved by Eggers and Hoppe \cite{EggersH}. It is natural to ask whether the prototype of generic singularity formation is of self-similar type.

First note that, by \eqref{23Jan12-J1}, for $s$ close to $0$, $T(s) \sim s^2$. Thus, to see the local picture of singularities, one should look at the scale $t \sim s^2$. In this scale, the Taylor expansion up to ``quartic terms'' of $\gamma$at $(0,0)$ is
\begin{align*}
\gamma(t,s)
	&= \beta(0)\,t + \beta'(0)\,ts + \frac{1}{2}\,\beta''(0)\,ts^2 + \frac{1}{6}\,\alpha'''(0)\,s^3\\
		&\qquad + \frac{1}{24}\,\alpha''''(0)\,s^4 + O(|t|^{5/2} + |s|^5).
\end{align*}
Now let
\[
e = \beta(0), u = \lr{\alpha'''(0)}{e^\perp} \text{ and } q = \lr{\beta'(0)}{e^\perp}.
\]
A simple computation leads to
\begin{align*}
\gamma(t,s)
	&= \big[t  - \frac{1}{2}\,q^2\,ts^2 - \frac{1}{8}\,qu\,s^4 \big]\,e
		 + \big[q\,ts + \frac{1}{6}\,u\,s^3\big]\,e^\perp\\
		&\qquad + O(|t|^{5/2} + |s|^5)\,e + O(|t|^2 + |s|^4)\,e^\perp.
\end{align*}
This shows that, infinitesimally, $\mcC_t$ is self-similar:
\[
\gamma(n^{-2}t, n^{-1}s) = \Big[n^{-2}t + n^{-4}\,\big(\lr{\gamma(t,s)}{e} - t\big)\Big]e
	+ n^{-3}\lr{\gamma(t,s)}{e^\perp}\,e^\perp
	+  O(n^{-5})\,e + O(n^{-4})\,e^\perp.
\]
%

\section{Timelike maximal surfaces in general vacuum spacetimes}

In this section we give the proof of Theorem \ref{Main2}. Let $(\mcM^{1 + 2},\stg)$ be a smooth oriented, time-oriented, globally hyperbolic Lorentzian manifold which satisfies the Einstein vacuum equation:
\begin{equation}
\stRic_{\alpha\beta} = 0.
	\label{VacBackground}
\end{equation}
Here $\stRic$ is the Ricci curvature of $\stg$. Let $t$ be a global time function on $\mcM$.

\subsection{Adapted coordinates}

Consider in $\mcM$ a closed spacelike acausal embedded curve $\mcC$ and a timelike (embedded) surface $\mcS$ which contains $\mcC$. We claim that $\mcC$ is a Cauchy curve for $\mcS$. Indeed, since $\mcC$ is acausal, it suffices to show that each inextendible causal curve in $\mcS$ must intersect $\mcC$. Let $\lambda$ a inextendible causal curve in $\mcS$. Since $\mcC$ is compact, the range of $t|_{\mcC}$ is bounded. Also, by the global hyperbolicity of $\mcM$, $t\big|_{\lambda}$ can attain any value in $\RR$, and in particular the value zero. The last two statement evidently imply that $\lambda$ intersects $\mcC$. The claim is proved.

By the above claim, $\mcS$ is globally hyperbolic (this can also be seen by noting that $t\big|_{\mcS}$ defines a time function on $\mcS$) and $\mcS$ is homeomorphic to $\RR \times \mcC$. 

Next, we follow Kulkarni \cite{Kulkarni} to define a `canonical' parametrization of $\mcS$ as follows. Assume that $\mcC$ is parametrized by $\{\gamma(s): s \in [0,\Xi]\}$ (where $\gamma(0) = \gamma(\Xi)$).  Let $\widehat\mcS \approx \RR \times \RR$ be the universal cover of $\mcS$. Let $\widehat\mcC$ be the lift of $\mcC$ and $\hat s$ be the lift of the parameter $s$. For any point $\hat p \in \widehat \mcS$, the null lines passing through $p$ intersect $\widehat\mcC$ at $p_-$ and $p_+$ whose $\hat s$-parameters are $x_-$ and $x_+$ where $x_- \leq x_+$. We then set
\[
\hat\tau(p) = \left\{\begin{array}{ll}
	0 &\text{ if } p \in \widehat\mcC,\\
	\frac{1}{2}(x_+ - x_-) &\text{ if $p$ is in the future of $\widehat\mcC$},\\
	-\frac{1}{2}(x_+ - x_-) &\text{ if $p$ is in the past of $\widehat\mcC$},
\end{array}\right.
\]
and
\[
\hat\xi(p) = \frac{1}{2}(x_+ + x_-).
\]
Then $(\hat\tau,\hat\xi)$ defines a global parametrization of $\widehat \mcS$. This descents to a parametrization $(\tau,\xi)$ of $\mcS$. Note that $\mcC = \{\tau = 0\}$.

Now, note that both $\partial_\tau + \partial_\xi$ and $\partial_\tau - \partial_\xi$ are null. Hence, the metric $g$ induced by $\stg$ on $\mcS$ takes the form
\[
g = A(-d\tau^2 + d\xi^2)
\]
where $A$ is nowhere zero. Since $\xi \equiv s$ on $\mcC$, which is spacelike, $A$ is positive. We thus have
\[
g = e^{2u(\tau,\xi)}\,[-d\tau^2 + d\xi^2].
\]
Here and in the rest of the paper, $\xi$ is assumed to take values in $\RR$ and all functions are periodic in $\xi$ with a fixed period $\Xi > 0$.

Near $\mcS$, we can complete $\{\tau,\xi\}$ to a local coordinate system $\{\rho,t,\xi\}$ such that $\mcS$ is at $\rho = 0$, $\partial_\rho$ is normal to $\mcS$ and $\stg(\partial_\rho,\partial_\rho) = 1$. We thus have
\begin{align}
\stg 
	&= - (e^{2u(t,\xi)} + 2\rho\,M(\rho,t,\xi))\,dt^2 + 4\rho\,N(\rho,t,\xi)\,dt\,d\xi + (e^{2u(t,\xi)} + 2\rho\,P(\rho,\tau,\xi))\,d\xi^2\nonumber\\
		&\qquad\qquad + d\rho^2 + 2\rho\,(Q(\rho,t,\xi)\,dt + S(\rho,t,\xi)\,d\xi).
	\label{BackgroundMetric}
\end{align}
Here all functions depending on $\rho$ are smooth up to $\rho = 0$.

\subsection{The governing equations}

We assume henceforth that $\mcS$ is maximal.

It is easy to see that the second fundamental form of $\mcS$ is 
\begin{equation}
h = -M(\tau,\xi)\,d\tau^2 + 2N(\tau,\xi)\,d\tau\,d\xi +  P(\tau,\xi)\,d\xi^2,
	\label{2ndForm}
\end{equation}
where, by a standard abuse of notations,
\[
M(\tau,\xi) = M(0,\tau,\xi), N(\tau,\xi) = N(0,\tau,\xi) \text{ and } P(\tau,\xi) = P(0,\tau,\xi).
\]

Since $\mcS$ is maximal, we thus have
\begin{equation}
H = \tr_{g} h = e^{-2u}\,M + e^{-2u}\,P = 0 \text{ along } \mcS.
	\label{ZeroMCurv}
\end{equation}

We next derive the constraint equations on $\mcS$. Let $\stnabla$ and $\nabla$ denote the Levi-Civita connection of $\stg$ and $g$, respectively, and $\stRiem$ denote the curvature tensor of $\stg$ on $\mcM$,
\[
\stRiem(X,Y,Z,W) = \stg((\stnabla_X \stnabla_Y - \stnabla_Y \stnabla_X - \stnabla_{[X,Y]})Z,W).
\]
By the Codazzi equation and the Einstein vacuum equation \eqref{VacBackground},
\begin{align*}
&0
	= \stRic(\partial_\xi,\partial_\rho)
	= -e^{2u}\stRiem(\partial_\tau, \partial_\xi,\partial_\rho, \partial_\tau) 
	= -e^{-2u}[\nabla_\tau h(\partial_\xi,\partial_\tau) - \nabla_\xi h(\partial_\tau, \partial_\tau)],\\
&0
	= \stRic(\partial_\tau,\partial_\rho)
	= e^{-2u}\,\stRiem(\partial_\xi,\partial_\tau,\partial_\rho, \partial_\xi) 
	= e^{-2u}[\nabla_\xi h(\partial_\tau,\partial_\xi) - \nabla_\tau h(\partial_\xi, \partial_\xi)].
\end{align*}
Rewriting using \eqref{ZeroMCurv}, we get
\begin{align}
0
	&= \nabla_\tau h(\partial_\xi,\partial_\tau) - \nabla_\xi h(\partial_\tau, \partial_\tau)\nonumber\\
	&= \partial_\tau N + \partial_\xi M + h_{\tau\tau}\,\stGamma_{\tau\xi}^\tau + h_{\tau\xi}(\stGamma_{\tau\xi}^\xi - \stGamma_{\tau\tau}^\tau) - h_{\xi\xi}\stGamma_{\tau\tau}^\xi\nonumber\\
	&= \partial_\tau N + \partial_\xi M - M\,\stGamma_{\tau\xi}^\tau + N(\stGamma_{\tau\xi}^\xi - \stGamma_{\tau\tau}^\tau) - P\stGamma_{\tau\tau}^\xi\nonumber\\
	&= \partial_\tau N + \partial_\xi M
	,\label{Codazzi1}\\
0
	&= \nabla_\xi h(\partial_\tau,\partial_\xi) - \nabla_\tau h(\partial_\xi, \partial_\xi)\nonumber\\
	&= \partial_\xi N - \partial_\tau P + h_{\tau\tau}(-\stGamma_{\xi\xi}^\tau) + h_{\tau\xi}(\stGamma_{\xi\tau}^\tau - \stGamma_{\xi\xi}^\xi) + h_{\xi\xi}\stGamma_{\xi\tau}^\xi\nonumber\\
	&= \partial_\xi N - \partial_\tau P + M\,\stGamma_{\xi\xi}^\tau + N(\stGamma_{\xi\tau}^\tau - \stGamma_{\xi\xi}^\xi) + P\,\stGamma_{\xi\tau}^\xi\nonumber\\
	&= \partial_\xi N - \partial_\tau P
	.\label{Codazzi2}
\end{align}

Next, by the Gauss equation and the Einstein vacuum equation \eqref{VacBackground},
\begin{align*}
0 
	&= -e^{-2u}\stRic(\partial_\tau,\partial_\tau) + e^{-2u}\stRic(\partial_\xi,\partial_\xi) -\stRic(\partial_\rho,\partial_\rho)\\
	&= -2\,e^{-4u}\stRiem(\partial_\xi, \partial_\tau, \partial_\tau, \partial_\xi)\\
	&= 2K + 2\,e^{-4u}[ h_{\xi\xi}\,h_{\tau\tau} - h_{\xi\tau}^2]\\
	&= 2K - 2\,e^{-4u}( M\,P + N^2).
\end{align*}
Here $K$ is the Gaussian curvature of $\mcS$,
\begin{align*}
K 
	&= -e^{-2u}[-\partial_{\tau\tau}u + \partial_{\xi\xi}u].
\end{align*}
We thus have, by \eqref{ZeroMCurv},
\begin{equation}
-\partial_{\tau\tau}u + \partial_{\xi\xi}u =  e^{-2u}(M^2 - N^2).
	\label{Gauss}
\end{equation}

To summarize, we have derived the following equations, which holds on $\mcS$,
\begin{align*}
&M + P = 0
	,\\
&\partial_\tau N + \partial_\xi M= 0
	,\\
&\partial_\xi N + \partial_\tau M = 0
	,\\
&-\partial_{\tau\tau}u + \partial_{\xi\xi}u =  e^{-2u}(M^2 - N^2)
	.
\end{align*}

\subsection{A blow up result}

We now give the proof of Theorem \ref{Main2}. 

First, note that $U = e^{-u}\,\partial_\xi$, $V = e^{-u}\,\partial_\tau$ and $\nu = \pm\partial_\rho$. By hypothesis,
\begin{align*}
0 
	< \stg(\stnabla_{e^{-u}\,\partial_\xi} (e^{-u}\,\partial_\xi),\nu)^2 - \stg(\stnabla_{e^{-u}\,\partial_\xi}(e^{-u}\,\partial_\tau),\nu)^2
	= e^{-4u}\,(P^2 - N^2) \text{ along } \mcC.
\end{align*}
Thus, \eqref{ZeroMCurv} implies that
\[
M^2 - N^2 > 0 \text{ on } \mcC.
\]
In particular, both $M - N$ and $M + N$ do not change sign on $\mcC$. On the other hand, by \eqref{Codazzi1}, \eqref{Codazzi2} and \eqref{ZeroMCurv},
\[
\partial_\tau(M + N) + \partial_\xi(M + N) = -\partial_\tau(M - N) + \partial_\xi(M - N) = 0,
\]
which implies that
\[
M(\tau,\xi) + N(\tau,\xi) = M(0,\xi - \tau) + N(0,\xi - \tau) \text{ and } M(\tau,\xi) - N(\tau,\xi) = M(0,\xi + \tau) + N(0,\xi + \tau).
\]
We conclude from the above discussion that both $M + N$ and $M - N$ do not change sign along $\mcS$ and are periodic in $\tau$. As $M^2 - N^2 > 0$ on $\mcC$, it follows that
\[
M^2 - N^2 > a > 0 \text{ along } \mcS.
\]
Recalling \eqref{Gauss}, we arrive at
\begin{equation}
-u_{,\tau\tau} + u_{,\xi\xi} \geq a\,e^{-2u}.
	\label{uBlowupEqn}
\end{equation}

We now follow a standard ODE technique to show that $u$ blows up in finite time (in either the future or the past or both). Let
\[
w(\tau) = \frac{1}{\Xi}\int_0^\Xi u(\tau,s)\,ds.
\]
Reversing time orientation if necessary, we can assume that $w'(0) \leq 0$.

We have
\[
w'' = \frac{1}{\Xi}\int_0^\Xi u_{\tau\tau}\,d\xi = \frac{1}{\Xi}\int_0^\Xi (u_{\tau\tau} - u_{\xi\xi})\,d\xi \leq -\frac{a}{\Xi}\int e^{-2u} \leq -a\,e^{-2w}.
\]
This implies that $w'$ is a strictly decreasing function. As $w'(0) \leq 0$, we thus have that $w'(\tau) < 0$ for all $\tau > 0$. We hence get
\[
\frac{d}{d\tau} (w')^2 \geq a\,\frac{d}{d\tau}e^{-2w}.
\]
In other words, $(w')^2 - a\,e^{-2w}$ is increasing. In particular,
\[
(w'(\tau))^2 \geq a\,e^{-2w(\tau)} - c_1 \text{ for all $\tau \geq 0$ and some constant $c_1 > 0$.}
\]

Using the differential inequality $w'' \leq a\,e^{-2w}$ and that $w'$ is strictly decreasing, we can find $\tau_0 > 0$ and $\delta > 0$ such that $w'(\tau) < -\delta$  for all $\tau > \tau_0$. This implies that for $\tau_1 > 0$ sufficiently large, 
\[
a\,e^{-2w(\tau)} > c_1 \text{ for } \tau \geq \tau_1.
\]
As $w' < 0$, the last two displayed inequalities give
\[
w'(\tau) \leq -\sqrt{a}\,e^{-w(\tau)} + c_2 < 0 \text{ for all $\tau \geq \tau_1$ and some constant $c_2 > 0$.}
\]
This implies that $\ln\frac{e^{-w}}{\sqrt{a}\,e^{-w} - c_2}$ is differentiable for $\tau \geq \tau_1$ and
\[
\frac{d}{d\tau}\left(\ln\frac{e^{-w}}{\sqrt{a}\,e^{-w} - c_2}\right) \leq -c_2 \text{ for all $\tau \geq \tau_1$}.
\]
Therefore, with $c_3 = \ln\frac{e^{-w(\tau_1)}}{\sqrt{a}\,e^{-w(\tau_1)} - c_2}$, there holds
\[
\ln\frac{e^{-w}}{\sqrt{a}\,e^{-w} - c_2} \leq c_3-c_2\tau \text{ for } \tau \geq \tau_1.
\]
As the left hand side is bounded from below by $-\frac{1}{2}\ln a$, this results in a contradiction for large $\tau$.

\appendix
\section{Regular timelike cylindrical maximal surfaces in $\RR^{1 + 3}$}

In this appendix, we construct a regular timelike cylindrical maximal surface in $\RR^{1 + 3}$. The arguments in Sections \ref{ssec:TSurfFlat} and \ref{ssec:TheEqns} remain valid in this context. Therefore, it suffices to point out a map 
\begin{align*}
\gamma: \RR \times \RR &\rightarrow \RR^3\\
	(t,s) &\mapsto \gamma(t,s)
\end{align*}
which is periodic in the $s$-factor such that
\begin{equation}
\left\{\begin{array}{l}
\Box\gamma = -\gamma_{,tt} + \gamma_{,ss} = 0,\\
\gamma_{,t}(0,s) = 0,\\
|\gamma_{,s}(0,s)| = 1,\\
|\gamma_{,s}(t,s)| \neq 0.
\end{array}\right.
	\label{10Oct11-Eq1}
\end{equation}
As before, let $\alpha(s) = \gamma(0,s)$. The third equation in \eqref{10Oct11-Eq1} says that $\alpha$ is parametrized by arclength. The fourth equation in \eqref{10Oct11-Eq1} guarantees that the map $(t,s) \mapsto (t,\gamma(t,s))$ is a regular immersion of $\RR \times \Sphere^1$ into $\RR^{1 + 3}$.

Given $\alpha$ which defines a parametrization by arclength of a closed curve in $\RR^3$, the first two equations in \eqref{10Oct11-Eq1} are solved by
\[
\gamma(t,s) = \frac{1}{2}[\alpha(s + t) + \alpha(s - t)].
\]
Thus the last equation in \eqref{10Oct11-Eq1} is equivalent to
\begin{equation}
\alpha_{,s}(s_1) + \alpha_{,s}(s_2) \neq 0 \text{ for all } s_1, s_2.
\label{10Oct11-Eq2}
\end{equation}
To finish the proof, we need to point out a curve $\alpha$ in $\RR^3$ which satisfies \eqref{10Oct11-Eq2}.

Let $P_1, P_2, P_3, P_4$ be four points in $\RR^3$ which don't belong to the same plane (but otherwise arbitrary). The curve $\alpha$ will be a smoothing of the piecewise linear curve $P_1P_2P_3P_4P_1$. The smoothing will be such that
\begin{itemize}
  \item[(LP)] The arc $C_i$ connecting $P_iP_{i+1}$ to $P_{i+1}P_{i+2}$ lies entirely in the plane $P_iP_{i+1}P_{i+2}$ for $i = 1, 2, 3, 4$. (Here we use the convention that $P_{5} = P_1$ and $P_6 = P_2$.)
  \item[(AC)] With respect to the plane $P_iP_{i+1}P_{i+2}$, the unit tangent to $C_i$ can be written as $e^{i\psi}$ where $\psi$ is a strictly monotone function whose image is contained in some interval $(a_i,a_i + \pi)$.
\end{itemize}
Evidently such curve $\alpha$ exists. A sketch of it follows.
\begin{figure}[h]
\begin{center}
\begin{tikzpicture}
\draw[fill=black] (0,0) circle (2pt)
(1,2) circle (2pt)
(3,1.5) circle (2pt)
(4,0) circle (2pt);
\draw[fill=black] (0,-0.4) node {$P_1$}
(1,2.4) node {$P_4$}
(3,1.9) node {$P_2$}
(4,-0.4) node {$P_3$};
\draw[dotted,thick] (0,0) -- (0.45,0.225)
(2.55,1.275) -- (3,1.5)
(3,1.5) -- (3.15,1.275)
(3.85,0.225) -- (4,0)
(4,0) -- (3.55,0.3)
(1.45,1.7) -- (1,2)
(1,2) -- (0.85,1.7)
(0.15,0.3) -- (0,0);
\draw (0.45,0.225) -- (2.22,1.11)
(2.34,1.17) -- (2.55,1.275)
(3.15,1.275) -- (3.85,0.225)
(3.55,0.3) -- (1.45,1.7)
(0.85,1.7) -- (0.15,0.3);
\draw (0.45,0.225) .. controls (0.225,0.1125) and (0.075,0.15) .. (0.15,0.3);
\draw (2.55,1.275) .. controls (2.775,1.3875) and (3.075,1.3875) .. (3.15,1.275);
\draw (3.85,0.225) .. controls (3.925,0.1125) and (3.775,0.15) .. (3.55,0.3);
\draw (1.45,1.7) .. controls (1.225,1.85) and (0.925,1.85) .. (0.85,1.7);
\end{tikzpicture}
\end{center}
\caption{An example of an initial curve which produces regular maximal surface.}
\end{figure}
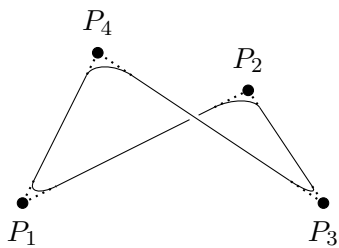

Note that \eqref{10Oct11-Eq2} is equivalent to the statement that $U(p) + U(q) \neq 0$ for any $p, q \in \alpha$ where $U$ is the unit tangent map of $\alpha$. Arguing by contradiction, assume that there is some $p$ and $q$ such that $U(p) + U(q) = 0$. If $p$ belongs to one of the straight segment connecting the $P_i$'s, say $P_1P_2$, then $U(p)$ is the same as $U(p')$ where $p'$ is the initial point of the junction arc $C_1$ connecting $P_1P_2$ to $P_2P_3$. We can thus assume without loss of generality that $p$ belongs to the arc $C_1$. It follows that $U(q)$ belongs to the plane $P_1P_2P_3$. By condition $(LP)$, it follows that $q$ must be either in segment $P_1P_2$, $P_2P_3$ or the arc $C_1$. As above, we can assume that $q$ belongs to the arc $C_1$. We thus have two points $p$ and $q$ on $C_1$ such that $U(p) + U(q) = 0$, but this contradicts condition $(AC)$.

The argument above shows that $\alpha$ satisfies \eqref{10Oct11-Eq2}. Hence, the timelike maximal surface containing $\alpha$ and perpendicular to the time slice $\RR^3$ containing $\alpha$ is a smooth (immersed) maximal surface in $\RR^{1 + 3}$.


\end{document}